\documentclass[10pt,a4paper]{article}
\usepackage[utf8]{inputenc}
\usepackage[english]{babel}
\usepackage{amsmath}
\usepackage{amsfonts}
\usepackage{amssymb}
\usepackage{amsthm}
\usepackage{float}
\usepackage{titlesec}

\usepackage[left=2cm,right=2cm,top=2cm,bottom=2cm]{geometry}

\numberwithin{equation}{section}

\newtheorem{lem}{Lemma}
\newtheorem{thm}{Theorem}
\newtheorem{cor}{Corollary}
\newtheorem{prop}{Proposition}

\theoremstyle{remark}
\newtheorem{rem}{\bf Remark}

\newcommand{\R}{\mathbb{R}}
\newcommand{\Sp}{\mathbb{S}}
\newcommand{\eqdef}{\overset{def}{=}}	

\newcommand{\tW}{\widetilde{W}}

\newcommand{\supp}{\mathrm{supp}}
\newcommand{\sgn}{\mathrm{sgn}}
\newcommand{\tG}{\widetilde{G}}
\newcommand{\tH}{\widetilde{H}}
\newcommand{\tf}{\widetilde{f}}

\title{A breakdown of injectivity for \\
weighted ray transforms in multidimensions}

\author{F.O. Goncharov\thanks{CMAP, Ecole Polytechnique, CNRS, Universit\'{e} Paris-Saclay, 91128, Palaiseau, France; \newline \indent email: fedor.goncharov.ol@gmail.com} 
	\and	 
		R. G. Novikov$^*$\footnote{IEPT RAS, 117997  
			Moscow, Russia;
			\newline\indent email: roman.novikov@polytechnique.edu}}
						
\begin{document}
\maketitle

\abstract{We consider weighted ray-transforms $P_W$ (weighted Radon transforms along straight lines) in $\R^d, \, d\geq 2,$ with strictly positive weights $W$. We construct an example of such a transform with non-trivial kernel in the space of infinitely smooth compactly supported functions on $\R^d$. 
In addition, the constructed weight $W$ is rotation-invariant continuous and is infinitely smooth almost everywhere on $\R^d\times \Sp^{d-1}$. In particular, by this construction we give counterexamples to some well-known injectivity results for weighted ray transforms for the case when the regularity of $W$ is slightly relaxed. We also give examples of continous strictly positive $W$ such that $\dim \ker P_W \geq n$ in the space of infinitely smooth compactly supported functions on $\R^d$ for arbitrary $n\in \mathbb{N}\cup \{\infty\}$, where $W$ are infinitely smooth for $d=2$ and infinitely smooth almost everywhere for $d\geq 3$.
\\


 }

\fontsize{12pt}{14pt}\selectfont

\section{Introduction}
\indent We consider the weighted ray transforms $P_W$ defined by 
\begin{align}\label{rad.def}
	&P_Wf(x,\theta) =  \int\limits_{\R}W(x + t\theta,\theta)f(x+ t\theta)
	\, dt, \, (x,\theta)\in T\Sp^{d-1},  \, d\geq 2,\\
	\label{ray_manifold}
	&T\Sp^{d-1} = \{(x,\theta)\in \R^d\times \Sp^{d-1} : x\theta = 0\},
\end{align}
where $f=f(x), \, W=W(x,\theta), \, x\in \R^d, \, \theta\in \Sp^{d-1}$.
Here, $W$ is the weight, $f$ is a test function on $\R^d$. In addition, we interpret $T\Sp^{d-1}$ as the set of all rays in $\R^d$. As a ray $\gamma$ we understand a straight line with fixed orientation. If $\gamma = \gamma(x,\theta), \, (x,\theta)\in T\Sp^{d-1}$, then
\begin{align}\label{ray_in_Rd}
	\begin{split}
		&\gamma(x,\theta) = \{y\in \R^d : y = x+t\theta, \, t\in \R\} 
		\text{ (up to orientation)},\\
		&\text{where }\theta \text{ gives the orientation of } \gamma.
	\end{split}
\end{align} 
\par We assume that 
\begin{align}\label{w.compl.}
	W = \overline{W} \geq c > 0, \, W\in L^{\infty}(\R^d\times \Sp^{d-1}),
\end{align}
where $\overline{W}$ denotes the complex conjugate of $W$, $c$ is a constant.
\par Note also that 
\begin{equation}\label{w_ray_form}
	P_Wf(x,\theta) = \int\limits_{\gamma} W(x,\gamma)f(x)\, dx, 
	\, \gamma = \gamma(x,\theta), 
\end{equation}
where 
\begin{equation}\label{w_gamma_def}
	W(x,\gamma) = W(x,\theta) \text{ for } x\in \gamma, \, \gamma = \gamma(x,\theta), \, 
	(x,\theta)\in T\Sp^{d-1}.
\end{equation}
\par The aforementioned transforms $P_W$ arise in various domains of pure and applied mathematics; see  \cite{lavrientiev1973first}, 
\cite{tretiak1980exponential}, 
\cite{quinto1983invertibility}, \cite{beylkin1984inversion}, 
\cite{markoe1985elementary}, \cite{davidfinch1986}, \cite{boman1987support}, 
\cite{sharafutdinov1992problem}, 
\cite{kunyansky1992generalized}, 
\cite{boman1993support}, \cite{boman1993example},
\cite{sharafutdinov1993uniqueness},   
\cite{kuchment1995local}, 
\cite{palamodov1996inversion},
\cite{arbuzov1998twodimensional},  
\cite{natterer2001mathematics}, 
\cite{novikovnonabel2002},
\cite{novikov2002inversion}, 
\cite{boman2004novikov}, 
\cite{bal2009inverese}, 
\cite{gindikin2010remark},
\cite{bal2011combined}, 
\cite{puro2013cormack}, 
\cite{novikov2014weighted}, \cite{ilmavirta2016coherent}, \cite{nguyen2017strength} and 
references therein.
\par In particular, the related results are the most developed for the case when $W\equiv 1$. In this case $P_W$ is reduced to the classical ray-transform $P$ (Radon transform along straight lines). The transform $P$ arises, in particular, in the X-ray transmission 
tomography. We refer to \cite{radon1917}, \cite{john1938ultrahyperbolic}, \cite{cormack1964representation}, 
\cite{gel1982integral}, \cite{helgason2001differential}, 
\cite{natterer2001mathematics} and references therein in connection with basic results for this classical case.
\par At present, many important results on transforms $P_W$ with other weights $W$ satisfying \eqref{w.compl.}  are also known; see the publications mentioned above with non-constant $W$ and references therein.
\par In particular, assuming \eqref{w.compl.} we have the following 
injectivity results.
\paragraph{\large Injectivity 1} \hspace*{-0.3cm}
(see \cite{davidfinch1986}). Suppose that $d\geq 3$ and $W\in C^2(\R^d\times \Sp^{d-1})$.
Then $P_W$ is injective on $L_0^p(\R^d)$ for $p > 2$, where $L_0^p$ denotes compactly
supported functions from $L^p$.
\paragraph{\large Injectivity 2}\hspace*{-0.3cm} (see \cite{markoe1985elementary}). Suppose that $d=2, \, W\in C^2(\R^2\times \Sp^1)$ and 
\begin{align}
	&0 < c_0 \leq W, \, \|W\|_{C^2(\R^2\times \Sp^1)} \leq N,
\end{align}
for some constants $c_0$ and $N$. Then, for any $p > 2$, there is $\delta = \delta(c_0, N, p) > 0$ such that $P_W$ is injective on $L^p(B(x_0,\delta))$ for any $x_0\in \R^2$, where
\begin{align}\nonumber
	\begin{split}
	&L^p(B(x_0,\delta)) = \{f\in L^p(\R^2) : \supp \, f \subset \overline{B}(x_0,\delta)\}, \\
	&\overline{B}(x_0,\delta) = \{x\in \R^2 : |x-x_0| \leq \delta\}.
	\end{split}
\end{align}
\paragraph{\large Injectivity 3} \hspace*{-0.3cm} (see \cite{quinto1983invertibility}). Suppose that $d=2, \, W\in C^1(\R^2\times \Sp^1)$ and $W$ 
is rotation invariant (see formula \eqref{w.rot.inv} below). Then
$P_W$ is injective on $L^p_0(\R^2)$ for $p\geq 2$.\\

\noindent In a similar way with \cite{quinto1983invertibility}, we say that $W$ is rotation invariant if  and only if
\begin{align}
	\begin{split}\label{w_rot_gamma_inv_def}
	&W(x,\gamma) \text{ is independent of the orientation of } \gamma, \\
	&W(x,\gamma) = W(Ax,A\gamma) \text{ for } x\in \gamma, \, \gamma\in T\Sp^{d-1}, \, 
	A\in O(d), 	
	\end{split}
\end{align}
where $T\Sp^{d-1}$ is defined in \eqref{ray_manifold}, $O(d)$ denotes the group of orthogonal transformations of $\R^d$.
\par Note also that property \eqref{w_rot_gamma_inv_def} can be rewritten in the form  
 \eqref{w.rot.inv}, \eqref{u.rot.inv.prop} or \eqref{w.new.sym}, \eqref{tu.rot.inv.prop}; see Section~\ref{sect.prelim}.

\paragraph{\large Injectivity 4} \hspace*{-0.3cm} (see \cite{boman1987support}). Suppose that $d=2$, 
$W$ is real-analytic on $\R^2\times \Sp^{1}$. Then $P_W$ is injective on $L_0^p(\R^2)$ for 
$p\geq 2$.\\

\par Injectivity 1 is a global injectivity for $d\geq 3$. Injectivity 2 is a local injectivity for $d=2$. Injectivity 3 is a global injectivity for $d=2$ for the rotation invariant case. 
Injectivity 4 is a global injectivity for $d = 2$ for the real-analytic case.
\par The results of Injectivity 1 and Injectivity 2 remain valid with $C^{\alpha}, \, \alpha > 1,$ in place of $C^2$ in the assumptions on $W$; see \cite{ilmavirta2016coherent}.
\par Injectivity 1 follows from Injectivity 2 in the framework of the layer-by-layer reconstruction approach. See \cite{davidfinch1986}, \cite{novikovnonabel2002}, \cite{ilmavirta2016coherent} and references therein in connection with the layer-by-layer 
reconstruction approach for weighted and non-abelian ray transforms in dimension $d\geq 3$.
\par The work \cite{boman1993example} gives a counterexample to Injectivity 4 for $P_W$ in $C_0^{\infty}(\R^2)$ for the case when the assumption that $W$ is real-analytic is relaxed to the assumption that $W$ is infinitely smooth, where $C_0^{\infty}$ denotes infinitely smooth compactly supported functions.
\par In somewhat similar way with \cite{boman1993example}, in the present work we obtain counterexamples to Injectivity 1, Injectivity 2 and Injectivity 3 for the case  when the regularity of $W$ is slightly relaxed. In particular, by these counterexamples we continue related studies of \cite{markoe1985elementary}, \cite{boman1993example} and \cite{gonnov2017wrt}.
\par More precisely, in the present work we construct $W$ and $f$ such that 
\begin{align}
	P_Wf \equiv 0 \text{ on } T\Sp^{d-1}, \, d\geq 2,
\end{align}
where $W$ satisfies \eqref{w.compl.}, $W$ is rotation-invariant (i.e., satisfies \eqref{w_rot_gamma_inv_def}),
\begin{align}
	\begin{split}\label{wintr_dthree} 
	&W \text{ is infinitely smooth almost everywhere on }\R^d\times \Sp^{d-1} 
	\text { and } \\
	&W\in C^{\alpha}(\R^d\times \Sp^{d-1}), 
	\text { at least, for any } \alpha\in (0,\alpha_0), \text{ where } \alpha_0= 1/16;
	\end{split}\\
	\label{fintr_dthree}
	\begin{split} 
	&f 
	\text { is a non-zero spherically symmetric infinitely 
	smooth and}\\
	&\text{compactly supported function on } \R^d;
	\end{split}
\end{align}
see Theorem~\ref{main.thm} of Section~\ref{sect.main}.

\par These $W$ and $f$ directly give the aforementioned  counterexamples to Injectivity 1 and Injectivity 3.
\par Our counterexample to Injectivity 1 is of particular interest (and is rather surprising) in view of the fact that the problem of finding $f$ on $\R^d$ from $P_Wf$ on $T\Sp^{d-1}$ 
for known $W$ is 
strongly overdetermined for $d\geq 3$. Indeed,
\begin{align*}
	&\dim \R^d = d, \, \dim T\Sp^{d-1} = 2d-2, \\
	& d < 2d-2 \text{ for } d\geq 3.
\end{align*}
\noindent This counterexample to Injectivity 1 is also rather surprising in view of the aforementioned 
layer-by-layer reconstruction approach in dimension $d\geq 3$.
\par Our counterexample to Injectivity 3 is considerably stronger than the preceeding
counterexample of \cite{markoe1985elementary}, where $W$ is not yet continuous and 
is not yet strictly positive (i.e., is not yet separated from zero by a positive constant).
\par Using our $W$ and $f$ of \eqref{wintr_dthree}, \eqref{fintr_dthree} for $d=3$ we also obtain the aforementioned counterexample to Injectivity 2; see Corollary~\ref{thm_delta_non_inj} of Section~\ref{sect.main}.
\par Finally, in the present work we also give examples of $W$ satisfying \eqref{w.compl.} such that 
$\dim \ker P_W\geq n$ in $C_0^{\infty}(\R^d)$ for arbitrary $n\in \mathbb{N}\cup \{\infty\}$, 
where $W\in C^{\infty}(\R^2\times \Sp^{1})$ for $d=2$ and $W$ satisfy \eqref{wintr_dthree} 
for $d\geq 3$; see Theorem~\ref{corol_dimension_infty} of Section~\ref{sect.main}. To our 
knowledge, examples of $W$ satisfying \eqref{w.compl.}, where $\dim \ker P_W\geq n$ (for example in $L_0^2(\R^d)$) were not yet given in the literature even for $n=1$ in dimension $d\geq 3$ 
and even for $n=2$ in dimension $d=2$.
\par In the present work we adopt and develop considerations of the famous work \cite{boman1993example} and of our very recent work \cite{gonnov2017wrt}.\\
\noindent In Section~\ref{sect.prelim} we give some preliminaries and notations.\\
\noindent Main results are presented in detail in Sections~\ref{sect.main}.\\
\noindent Related proofs are given in Sections~\ref{sect_thm2_proof}-\ref{pr.lm5}.

\section{Some preliminaries}\label{sect.prelim}
\paragraph{\large Notations.} Let 
\begin{align}
\label{pr.ts2}
&\Omega = \R^d\times \Sp^{d-1},\\ \label{distance.line.func}
&r(x,\theta) = |x-(x\theta)\theta|, \, (x,\theta)\in \Omega, \\
 \label{pr.l0}
&\Omega_0(\delta) = \{(x,\theta) \in \Omega : r(x,\theta) > \delta \},\\ \label{pr.l1}
&\Omega_1(\delta) = \Omega\backslash \Omega_0(\delta) 
= \{(x,\theta) \in \Omega : r(x,\theta) \leq \delta\}, 
\, \delta > 0,\\ \label{loc.weight.domain}
&\Omega(\Lambda) = \{(x,\theta)\in \R^d\times \Sp^{d-1} : r(x,\theta)\in 
\Lambda\}, \, \Lambda \subset [0, +\infty),\\ \label{t0_def}
&T_0(\delta) = \{(x,\theta)\in T\Sp^{d-1} : |x| > \delta\},\\ \label{t1_def}
&T_1(\delta) = \{(x,\theta)\in T\Sp^{d-1} : |x| \leq \delta\}, \, \delta > 0,\\ 
\label{t_def}
&T(\Lambda) = \{(x,\theta)\in T\Sp^{d-1} : |x|\in \Lambda\}, \, 
\Lambda \subset [0, +\infty), \\
\label{int.def}
&\mathcal{J}_{r,\varepsilon}   = (r-\varepsilon, r+\varepsilon)\cap [0, +\infty), 
\, r\in [0, +\infty), \, \varepsilon > 0.
\end{align}
\par The set $T_0(\delta)$ in \eqref{t0_def} is considered as the set of all rays in $\R^d$ which are located at distance greater than $\delta$ from the origin.
\par The set $T_1(\delta)$ in \eqref{t1_def} is considered as the set of all rays in $\R^d$ which are located at distance less or equal than $\delta$.

\par We also consider the projection 
\begin{align}\label{pi_projection_1}
	&\pi : \Omega \rightarrow T\Sp^{d-1},\\
	&\pi(x,\theta) = (\pi_\theta x, \theta), \, (x,\theta)\in \Omega, \\ 
	\label{pi_projection_3}
	&\pi_\theta x = x - (x\theta)\theta.
\end{align}
\par In addition, $r(x,\theta)$ of \eqref{distance.line.func} is the distance from the origin $\{0\}\in \R^d$ to the ray $\gamma = \gamma(\pi(x,\theta))$ (i.e., $r(x,\theta) = |\pi_\theta x|$). The rays will be also denoted by 
\begin{equation}\label{ray_def_omega}
\gamma = \gamma(x,\theta) \eqdef \gamma(\pi(x,\theta)), \, (x,\theta)\in \Omega.
\end{equation}
\par We also consider 
\begin{equation}
	P_Wf(x,\theta) = P_Wf(\pi(x,\theta)) \text{ for } (x,\theta)\in \Omega.
\end{equation}
\par We also define
\begin{align}
\begin{split}
&B(x_0,\delta) = \{x\in \R^d : |x-x_0| < \delta\}, \\
&\overline{B}(x_0,\delta) = \{x\in \R^d : |x-x_0|\leq \delta\},\,x_0\in \R^d, \,  \delta > 0,\\ 
\end{split}\\
\label{ball_def}
&B = B(0, 1), \, \overline{B} = \overline{B}(0,1).
\end{align}
\par For a function $f$ on $\R^d$ we denote its restriction to a subset $\Sigma\subset
\R^d$ by $f|_\Sigma$.
\par By $C_0, \, C_0^{\infty}$ we denote continuous compactly supported and infinitely smooth compactly supported functions, respectively.
\par By $C^\alpha(Y), \,  \alpha\in(0,1)$, we denote the space of $\alpha$-H\"{o}lder functions 
on $Y$ with the norm:
\begin{align}
	\begin{split}
	&\|u\|_{C^\alpha(Y)} = \|u\|_{C(Y)} + \|u\|'_{C^\alpha(Y)},\\
	&\|u\|'_{C^\alpha(Y)} = \sup\limits_{\substack{y_1,y_2\in Y \\ 
												   |y_1-y_2| \leq 1}}
							\dfrac{|u(y_1)-u(y_2)|}{|y_1-y_2|^{\alpha}},	
	\end{split}
\end{align} 
where $\|u\|_{C(Y)}$ denotes the maximum of $|u|$ on $Y$.

\paragraph{\large Rotation invariancy.} 
Using formula \eqref{w_gamma_def}, for positive and continous $W$, property \eqref{w_rot_gamma_inv_def} can be rewritten in the following equivalent form:
\begin{equation}\label{w.rot.inv}
 W(x,\theta) = U(|x-(x\theta)\theta|, x\theta), \, x\in \R^d, \, \theta\in \Sp^{d-1},
\end{equation}
for some positive and continuous $U$ such that
\begin{equation}\label{u.rot.inv.prop}
  U(r,s) = U(-r,s) = U(r,-s), \, r\in \R, \, s\in \R.
\end{equation}
\noindent  In addition, symmetries \eqref{w.rot.inv}, \eqref{u.rot.inv.prop} of $W$ can be also written as
\begin{align}\label{w.new.sym}
	&W(x,\theta) = \tilde{U}(|x|,x\theta), \, (x,\theta)\in \Omega,\\ \label{tu.rot.inv.prop}
	&\tilde{U}(r,s) = \tilde{U}(-r,s) = \tilde{U}(r,-s), \, r\in \R, \, s\in \R.
\end{align}
where $\tilde{U}$ is positive and continuous on $\R\times \R$. Using the formula $|x|^2 = |x\theta|^2 + r^2(x,\theta),$ one can see that symmetries \eqref{w.rot.inv}, 
\eqref{u.rot.inv.prop} and symmetries \eqref{w.new.sym}, \eqref{tu.rot.inv.prop} of $W$ are equivalent.

\paragraph{\large Partition of unity.}
We recall the following classical result (see, e.g., Theorem 5.6 in \cite{do1992riemannian}):\\ 
{\itshape
	Let $\mathcal{M}$ be a $C^{\infty}$-manifold, which is Hausdorff and has a countable base. Let also $\{U_i\}_{i=1}^{\infty}$ be an open locally-finite cover of $\mathcal{M}$.
	\par Then there exists a $C^{\infty}$-smooth locally-finite partition of unity $\{\psi_i\}_{i=1}^{\infty}$ on $\mathcal{M}$,
	such that 
	\begin{equation}\label{gen.part.prop}
		\supp\, \psi_i\subset U_i.
	\end{equation}}
\vspace{-0.4cm}
\par In particular, any open interval $(a,b)\subset \R$ and $\Omega$ satisfy the conditions for $\mathcal{M}$ of this statement. It will be used in Subsection~\ref{subsect_construct_dtwo}.

\section{Main results}\label{sect.main}
\begin{thm}\label{main.thm}
	There exist a weight $W$ satisfying \eqref{w.compl.} and 
	a non-zero function \newline ${f\in C_0^\infty(\R^d), \, d\geq 2}$, such that 
	\begin{equation}\label{m.thm.stm.}
		P_Wf\equiv 0 \text{ on } T\Sp^{d-1},
	\end{equation}
	where $P_W$ is defined in \eqref{rad.def}. In addition, $W$ is rotation invariant, i.e., satisfies \eqref{w.rot.inv}, and $f$ is spherically symmetric with $\supp \, f\subseteq \overline{B}$.
	 Moreover, 
		\begin{align}\label{thm_prop_smooth}
			&W\in C^{\infty}(\Omega \backslash \Omega(1)),\\ \label{thm_prop_holder}
			&W \in C^\alpha(\R^d\times \Sp^{d-1}) \text{ for any } \alpha\in (0,\alpha_0),\, \alpha_0 = 1/16,\\ \label{thm_prop_pos}
			&W\geq 1/2 \text{ on }\Omega  \text{ and }W\equiv 1 \text{ on } \Omega([1,+\infty)),\\
			\label{thm_w_unity}
			&W(x,\theta) \equiv 1 \text{ for } |x| \geq R > 1, 	\, \theta\in \Sp^{d-1},
		\end{align}
	 where $\Omega, \, \Omega(1), \, \Omega([1,+\infty))$ are defined by  
	 \eqref{pr.ts2}, \eqref{loc.weight.domain}, $R$ is a constant.
\end{thm}
\par The construction of $W$ and $f$ proving Theorem~\ref{main.thm} is presented below in Subsections~\ref{subsect_construct_dtwo},~\ref{subsect.extension}. 
In addition, this construction consists of its version in dimension $d=2$ (see Subsection~\ref{subsect_construct_dtwo}) and its subsequent extension to the case of $d\geq 3$ (see Subsection~\ref{subsect.extension}). 
\par Theorem~\ref{main.thm} directly gives counterexamples to Injectivity 1 and Injectivity 3 
of Introduction. Theorem~\ref{main.thm} also implies the following counterexample to Injectivity 2 of Introduction:
\begin{cor}\label{thm_delta_non_inj}
For any $\alpha\in (0,1/16)$ there is $N > 0$ such that for any $\delta > 0$ 
there are $W_\delta, \, f_\delta$ satisfying
\begin{align}\label{cexample_inj2_start}
	&W_\delta \geq 1/2, \,W_\delta\in C^{\alpha}(\R^2\times \Sp^1), \,  \|W_\delta\|_{C^{\alpha}(\R^2\times \Sp^1)} \leq N\\
	\label{cexample_inj2_fsupp}
	&f_{\delta}\in C^{\infty}(\R^2), \, f_{\delta}\not\equiv 0, \, \supp\, f_{\delta}
	\subseteq \overline{B}(0,\delta),\\ \label{cexample_inj2_end} 
	&P_{W_\delta}f_{\delta} \equiv 0 \text{ on }T\Sp^1.
\end{align}
\end{cor}
\noindent The construction of $W_{\delta}, \, f_{\delta}$ proving Corollary~\ref{thm_delta_non_inj} is presented in  Subsection~\ref{proof_cor1}.
\begin{thm}\label{corol_dimension_infty}
	For any $n\in \mathbb{N}\cup \{\infty\}$ there exists a weight $W_n$ satisfying \eqref{w.compl.} such that
	\begin{equation}\label{dim_ker_w}
		\dim \ker P_{W_{n}}\geq n \text{ in } C_0^{\infty}(\R^d), \, d\geq 2,
	\end{equation}
	where $P_W$ is defined in \eqref{rad.def}.
	Moreover, 
	\begin{align}
    \label{wn_smooth_d2}
    &W_n \in C^{\infty}(\R^2\times \Sp^1) \text{ for }d=2, \\
    \label{wn_smooth_d3}
    \begin{split}
    	&W_n \text{ is infinitely smooth almost everywhere on }
    	\R^d\times \Sp^{d-1} \text{ and } \\
    &W_n\in C^{\alpha}(\R^d\times \Sp^{d-1}),\, \alpha\in (0,1/16) \text{ for } d\geq 3,
    \end{split}\\
    \label{wn_identity_outside}
	&W_n(x,\theta) \equiv 1 \text{ for } |x| \geq R > 1, \, 
	\theta\in \Sp^{d-1} \text{ for } n\in \mathbb{N}, \, d\geq 2,
    \end{align}
    where $R$ is a constant.
\end{thm}
\noindent The construction of $W_{n}$ proving Theorem~\ref{corol_dimension_infty} is 
presented in Section~\ref{sect_thm2_proof}. In this construction we proceed from Theorem~\ref{main.thm} of the present work for $d\geq 3$ and from the result of \cite{boman1993example} for $d=2$. In addition, for this construction it is essential that $n < +\infty$ in \eqref{wn_identity_outside}.

\subsection{Construction of $f$ and $W$ for $d=2$}\label{subsect_construct_dtwo}
\par In dimension $d=2$, the construction of $f$ and $W$ adopts and develops considerations of \cite{boman1993example} and \cite{gonnov2017wrt}. In particular, we construct $f$, first, and then $W$ (in this construction we use notations of Section~\ref{sect.prelim} for $d=2$). 
In addition, this construction is commented in Remarks~\ref{rem_constr_f}-\ref{rem_constr_w_local} below.
\paragraph*{\large Construction of $f$.}
The function $f$ is constructed as follows:
\begin{align}\label{m.def.ser.}
	&f = \sum\limits_{k=1}^{\infty}\dfrac{f_k}{k!},\\ \label{m.def.fk}
	&f_k(x) = \tf_k(|x|) = \Phi(2^k(1-|x|))\cos(8^k|x|^2), x\in \R^2, \, k\in \mathbb{N},
\end{align}
for arbitrary $\Phi\in C^{\infty}(\R)$ such that
\begin{align}\label{phi.suppint}
	&\supp\, \Phi = [4/5, 6/5], \\
	\label{phi.supp}
	&0 < \Phi(t) \leq 1 \text{ for } t\in (4/5,6/5), \\
	\label{phi.unity.cond}
	&\Phi(t) = 1, \text{ for } t\in [9/10, 11/10], \\
	\label{phi.monotonicity.cond}
	\begin{split}
		&\Phi \text{ monotonously increases on } [4/5,9/10]\\
		&\text{and monotonously decreases on } [11/10, 6/5].
	\end{split}
\end{align}

\noindent Properties \eqref{phi.suppint}, \eqref{phi.supp} imply that functions $\tf_k$ (and functions $f_k$) in \eqref{m.def.fk} have disjoint supports: 
\begin{align}\label{f_k_supports}
	\begin{split}
	&\supp \tf_i \cap \supp \tf_j = \emptyset \text{ if } i\neq j,\\
	&\supp \tf_k = [1-2^{-k}\left(\frac{6}{5}\right), 1-2^{-k}
	\left(\frac{4}{5}\right)],\, i,\, j,\, k\in\mathbb{N}.
	\end{split}
\end{align}
This implies the convergence of series in \eqref{m.def.ser.} for every fixed $x\in \R^2$. 
\begin{lem}\label{b.lem.ch.sgn}
	Let $f$ be defined by \eqref{m.def.ser.}-\eqref{phi.unity.cond}. Then $f$ is spherically symmetric, \newline $f\in C_0^{\infty}(\R^2)$ and $\supp \, f \subseteq \overline{B}$. In addition, if $\gamma\in T\Sp^{1}, \, \gamma\cap B\neq \emptyset$, then $f|_\gamma\not\equiv 0$ and $f|_\gamma$ has non-constant sign.
\end{lem}
\par Lemma~\ref{b.lem.ch.sgn} is similair to Lemma 1 of \cite{gonnov2017wrt} and it 
is, actually, proved in Section~4.1 of \cite{gonnov2017wrt}.

\begin{rem}\label{rem_constr_f}
	Formulas \eqref{m.def.ser.}-\eqref{phi.unity.cond} for $f$ are  
	similar to the formulas for $f$ in \cite{boman1993example}, where $P_W$ was considered in $\R^2$, and also to the formulas for $f$ in \cite{gonnov2017wrt}, 
	where the weighted Radon transform $R_W$ along hyperplanes was considered in $\R^3$.
	The only difference between \eqref{m.def.ser.}-\eqref{phi.unity.cond} and the related formulas in \cite{gonnov2017wrt} is the dimension $d=2$ in \eqref{m.def.ser.}-\eqref{phi.unity.cond} instead of $d=3$ in \cite{gonnov2017wrt}. At the same time, the important difference between \eqref{m.def.ser.}-\eqref{phi.unity.cond} and the related formulas in \cite{boman1993example} is that in formula \eqref{m.def.fk} the factor $\cos(8^k|x|^2)$ depends only on $|x|$, whereas in \cite{boman1993example} the corresponding factor is $\cos(3^k\phi)$ which depends only on the angle $\phi$ in the polar coordinates in $\R^2$.
In a similar way with \cite{boman1993example}, \cite{gonnov2017wrt}, we use the property that the restriction of the 
function $\cos(8^k|x|^2)$ to an arbitrary ray $\gamma$ intersecting the open ball oscillates sufficiently fast (with change of the sign) for large $k$.	
\end{rem}

\paragraph*{\large Construction of $W$.}
\par In our example $W$ is of the following form:
\begin{align}\label{b.W.form}
\begin{split}
	W(x,\theta) &= \phi_1(x) \left(
	\sum\limits_{i=0}^{N}\xi_i(r(x,\theta))W_i(x,\theta)
	\right) + \phi_2(x) \\
	& = 
	\phi_1(x)\left( \xi_0 (r(x,\theta)) W_0(x,\theta) + \sum\limits_{i=1}^{N}
	\xi_i(r(x,\theta))
	W_i(x,\theta)
	\right) + \phi_2(x)
	,\, (x,\theta)\in \Omega, 
\end{split}
\end{align}
where
\begin{align}
\begin{split}\label{part.unity.phi.def}
&\phi_1 = \phi_1(|x|), \phi_2 = \phi_2(|x|) \text{ is a $C^{\infty}$-smooth partition of unity on } \R^2 \text{ such that}, \\
&\phi_1 \equiv 0 \text{ for } |x|\geq R > 1, \, \phi_1 \equiv 1 \text{ for } |x|\leq 1,\\
&\phi_2 \equiv 0 \text{ for } |x|\leq 1,
\end{split}\\
\label{b.part.unit}
&\{\xi_i(s), \, s\in \R\}_{i=0}^{N} \text{ is a $C^{\infty}$-\, smooth partition of unity on }\R,\\
\label{b.un.part.prop}
&\xi_i(s) = \xi_i(-s), \, s\in \R, \, i=\overline{0,N},
\end{align}
\vspace{-0.5cm}
\begin{align}
\label{b.weight.desc}
\begin{split}
&W_i(x,\theta) \, \text{are bounded, continuous, strictly positive}\\
&\text{and rotation invariant (according to 
\eqref{w.rot.inv}), \eqref{tu.rot.inv.prop} on }\\ 
&\text{the open vicinities of } \supp\, \xi_i(r(x,\theta)), \, i = \, \overline{0,N}, \text{ respectively}.
\end{split}\hspace{2.0cm}
\end{align}
From the result of Lemma~\ref{b.lem.ch.sgn} and from \eqref{part.unity.phi.def} it follows that 
\begin{align}\label{ray.tr.of.form}
	\begin{split}
	P_Wf(x,\theta) &= \xi_0(|x|) P_{W_0}f(x,\theta) + 
	\sum\limits_{i=1}^N \xi_i(|x|) P_{W_i}f(x,\theta), \, (x,\theta)\in T\Sp^1, 
	\end{split}
\end{align}
where $W$ is given by \eqref{b.W.form}.

\par From \eqref{b.W.form}-\eqref{b.weight.desc} it follows that $W$ of \eqref{b.W.form} satisfies the conditions \eqref{w.compl.}, \eqref{w.new.sym}, \eqref{tu.rot.inv.prop}.

\par The weight $W_0$ is constructed in next paragraph and has the following properties:
\begin{align}\label{w0.c.rinv}
	&W_0\text{ is bounded, continuous and rotation invariant on } \Omega(1/2, +\infty),\\
	\label{w0_holder_property}
	\begin{split}
	&W_0\in C^{\infty}(\Omega\left((1/2, 1)\cup (1, +\infty)\right)) \text{ and }\\
	&W_0\in C^{\alpha}(\Omega(1/2, + \infty)) \text{ for }\alpha\in (0,1/16),
	\end{split}\\
	\label{w0.pos.d}
	\begin{split}
		&\text{there exists $\delta_0\in (1/2,1)$ such that: } \\
		&\qquad \qquad W_0(x,\theta) \geq 1/2 \text{ if } r(x,\theta) > \delta_0,\\
		&\qquad \qquad W_0(x,\theta) = 1 \quad \text{ if } r(x,\theta)\geq 1,
	\end{split}\\ \label{w0.zero}
	&P_{W_0}f(x,\theta) = 0 \text{ on } \Omega((1/2, +\infty)),
\end{align}
where $P_{W_0}$ is defined according to \eqref{rad.def} for $W = W_0$, $f$ is given by 
\eqref{m.def.ser.}, \eqref{m.def.fk}.
\par In addition, 
\begin{align}\label{supp.xi0}
	&\supp \, \xi_0 \subset (-\infty, -\delta_0)\cup (\delta_0, +\infty),\\
	\label{super.prop.xi0}
	&\xi_0(s) = 1 \text{ for } |s|\geq 1,
\end{align}
where $\delta_0$ is the number of \eqref{w0.pos.d}.
\par In particular, from \eqref{w0.pos.d}, \eqref{supp.xi0} it follows that
\begin{equation}\label{w0.positiv}
	W_0(x,\theta)\xi_0(r(x,\theta)) > 0 \text{ if } \xi_0(r(x,\theta)) > 0.
\end{equation}
\par In addition,
\begin{align}\label{w.loc.prop}
	& \xi_i(r(x,\theta))W_i(x,\theta) \text{ are bounded, rotation invariant and }C^{\infty} \text{ on }	\Omega,\\ 
	& W_i(x,\theta)\geq 1/2 \text{ if }  \xi_i(r(x,\theta)) \neq 0,\\
	\label{w.loc.prop.last}
	&P_{W_i}f(x,\theta) = 0 \text{ on } (x,\theta)\in T\Sp^1, \text{ such that } \xi_i(r(x,\theta))\neq 0,\\ \nonumber
	&i = \overline{1,N}, \, (x, \theta)\in\Omega.
\end{align}
Weights $W_1, \dots, W_N$ of \eqref{b.W.form} and $\{\xi_i\}_{i=0}^{N}$ are constructed in Subsection~\ref{subsect_construct_dtwo}.
\par Theorem~\ref{main.thm} for $d=2$ follows from Lemma~\ref{b.lem.ch.sgn} and formulas \eqref{b.W.form}-\eqref{w0.zero}, \eqref{w0.positiv}-\eqref{w.loc.prop.last}. 
\par We point out that the construction of $W_0$ of \eqref{b.W.form} is substantially different from the construction of $W_1, \dots , W_N$. 
The weight $W_0$ is defined for the rays $\gamma\in T\Sp^1$ which can be close to the boundary $\partial B$ of $B$ which results in  restrictions on global smoothness of $W_0$.

\begin{rem}\label{rem_general_constr_w}
 The construction of $W$ summarized above in formulas \eqref{b.W.form}-\eqref{w.loc.prop.last} arises in the framework of finding $W$ such that 
 \begin{align}\label{rw_zero_remark}
 	\begin{split}
 	&P_Wf \equiv 0 \text{ on }T\Sp^{1} \text{ for } f \text{ defined in 
 	\eqref{m.def.ser.}-\eqref{phi.monotonicity.cond}} ,
 	\end{split}
 \end{align}
 under the condition that $W$ is strictly positive, sufficiently regular and rotation invariant (see formulas \eqref{w.compl.}, \eqref{w.rot.inv}, \eqref{u.rot.inv.prop}).
 In addition, the weights $W_{i}, \, i = 0, \dots , N,$ in \eqref{b.W.form} are constructed in a such a way that
 \begin{equation}\label{rem_local_zeros}
 	P_{W_i}f = 0 \text{ on } V_i, \, i = 0, \dots, N, 
\end{equation}
under the condition that $W_{i} = W_i(x,\gamma)$ are strictly positive, sufficiently regular and rotation invariant for $x\in \gamma, \, \gamma\in V_i\subset T\Sp^1,\, i = 0, \dots, N$,   
where 
\begin{align}
	&\{V_i\}_{i=0}^{N} \text{ is an open cover of } T\Sp^{1} \text{ and } 
	V_0 = T_0(\delta_0),\\
	\label{rem_lambda_mention}
	&V_i = T(\Lambda_i) \text{ for some open } \Lambda_i\subset \R, \, 
	i = 0, \dots, N,
\end{align}
where $T_0$ is defined in \eqref{t0_def}, $\delta_0$ is the number of \eqref{w0.pos.d}, 
$T(\Lambda)$ is defined in \eqref{t_def}. In addition, the functions $\xi_i, \, i=0,\dots, N,$ in \eqref{b.W.form} can be interpreted as a partition of unity on $T\Sp^{1}$ subordinated to the open cover $\{V_i\}_{i=0}^{N}$. The aforementioned construction of $W$ is a two-dimensional analog of the construction developed in \cite{gonnov2017wrt}, where   the weighted Radon transform $R_W$ along hyperplanes was considered in $\R^3$. At the same time, 
the construction of $W$ of the present work is similar to the construction in \cite{boman1993example} with the important difference that in the present work $f$ is spherically symmetric and $W, \, W_i, \, i = 0, \dots, N$, are rotation invariant.
\end{rem}

\paragraph*{\large Construction of $W_0$.}

Let $\{\psi_k\}_{k=1}^{\infty}$ be a $C^{\infty}$ partition of unity on  
	$(1/2, 1)$ such that
\begin{align}\label{psik.unit.part.def}
	\begin{split}
	&\supp\, \psi_k\subset (1-2^{-k+1}, 1-2^{-k-1}), \, k\in \mathbb{N},
	\end{split}\\
	\label{psik.unit.part.deriv}
	&\text{first derivatives $\psi_k'$ satisfy the bounds: } \sup |\psi_k'| \leq C2^{k},	
\end{align}
where $C$ is a positive constant. Actually, functions $\{\psi_k\}_{k=1}^{\infty}$ satisfying \eqref{psik.unit.part.def}, \eqref{psik.unit.part.deriv} were used in  considerations of \cite{boman1993example}.
\par Note that
\begin{equation}
	1-2^{-(k-2)-1} < 1 -2^{-k}(6/5), \, k\geq 3.
\end{equation}
Therefore,
\begin{equation}\label{supp.coincide}
	\text{for all } s_0,\,t_0\in \R , s_0\in \supp\, \psi_{k-2}, \, t_0\in \supp\, \Phi(2^k(1-t))
	\Rightarrow s_0 < t_0, \, k\geq 3.
\end{equation}
\par Weight $W_0$ is defined  by the following formulas
\begin{align}\label{b.w0.def}
	&W_0(x,\theta) = 
	\begin{cases}
		1- G(x,\theta)\sum\limits_{k=3}^{\infty}k!f_k(x)
		\dfrac{\psi_{k-2}(r(x,\theta))}{H_k(x,\theta)}, \, 1/2 < r(x,\theta) < 1,\\
		1, \, r(x,\theta) \geq 1
	\end{cases}, \\ \label{b.gkhk.def}
	&G(x,\theta) = \hspace*{-0.2cm}\int\limits_{\gamma(x,\theta)}
	\hspace*{-0.3cm}f(y)\, dy, \, 
	H_k(x,\theta) = \hspace*{-0.2cm}\int\limits_{\gamma(x,\theta)}
	\hspace*{-0.3cm}f_k^2(y)\, dy, \, x\in \R^{2},\, \theta\in \Sp^{1},
\end{align}
where $f, \, f_k$ are defined in \eqref{m.def.ser.}, \eqref{m.def.fk}, respectively, 
rays $\gamma(x,\theta)$ are given by \eqref{ray_def_omega}. 
\par Formula \eqref{b.w0.def} implies that $W_0$ is defined on $\Omega_0(1/2)\subset \Omega$. \par Due to \eqref{m.def.fk}-\eqref{phi.unity.cond}, \eqref{psik.unit.part.def}, \eqref{supp.coincide}, in \eqref{b.gkhk.def} we have that 
\begin{align}\label{H_k_nonzero}
&H_k(x,\theta) \neq 0 \text{ if } \psi_{k-2}(r(x,\theta))\neq 0, \, 
(x,\theta)\in \Omega,\\
\label{H_k_smoothness}
&\dfrac{\psi_{k-2}(r(x,\theta))}{H_k(x,\theta)}\in C^{\infty}(\Omega(1/2, 1)),
\end{align}
where $r(x,\theta)$ is defined in \eqref{distance.line.func}, 
$\Omega, \, \Omega(\cdot)$ are defined in \eqref{pr.ts2}, \eqref{loc.weight.domain}, $d=2$.
\par Also, for any fixed $(x,\theta)\in \Omega, \, 1/2<r(x,\theta)$, the series in the right hand-side of \eqref{b.w0.def} has only a finite number of non-zero terms (in fact, no more than two) and, hence, the weight $W_0$ is well-defined.
\par By the spherical symmetry of $f$, functions $G, H_k$ in \eqref{b.w0.def} are of the type \eqref{w.rot.inv} (and \eqref{w.new.sym}). Therefore, $W_0$ is rotation invariant (in the sense of \eqref{w.rot.inv} and \eqref{w.new.sym}). 
\par Actually, formula \eqref{w0.zero} follows from \eqref{m.def.ser.}, \eqref{m.def.fk}, \eqref{b.w0.def}, \eqref{b.gkhk.def} (see Subsection~\ref{subs.w0.zero} for details).
\par Using the construction of $W_0$ and the assumption that $r(x,\theta) > 1/2$ one can see that $W_0$ is $C^{\infty}$ on its domain of definition, possibly, except points with $r(x,\theta)=1$.

\par Note also that due to \eqref{m.def.ser.}, \eqref{m.def.fk}, the functions $f_k, \, G, H_k$,
used in \eqref{b.w0.def}, \eqref{b.gkhk.def} can be considered as functions of one-dimensional arguments.
\par Formulas \eqref{w0.c.rinv}-\eqref{w0.pos.d} are proved in Subsection~\ref{subsect_pr_lem2}.

\begin{rem}\label{rem_constr_w0}
\par Formulas \eqref{b.w0.def}, \eqref{b.gkhk.def} given above for the weight $W_0$ are considered for the rays from $T_0(\delta_0)$ (mentioned in Remark~\ref{rem_general_constr_w}) and, in particular, for rays close to the tangent rays to $\partial B$. 
These formulas are direct two-dimensional analogs of the related formulas in \cite{gonnov2017wrt}. At the same time, formulas \eqref{b.w0.def}, \eqref{b.gkhk.def} are similar to the related formulas in \cite{boman1993example} with the important difference  that $f, \, f_k$ are spherically symmetric in the present work and, as a corollary, $W_0$ is rotation invariant. Also, in a similar way with \cite{boman1993example}, \cite{gonnov2017wrt}, in the present work we show that $G(x,\theta)$ tends to zero sufficiently fast as $r(x,\theta)\rightarrow 1$. This is a very essential point for continuity of $W_0$ and it is given in Lemma~\ref{b.cont.lem} of Subsection~\ref{subsect_pr_lem2}.
\end{rem}

\paragraph*{\large Construction of $W_1, \dots, W_N$ and $\xi_0,\dots, \xi_N$}
\begin{lem}\label{b.lem.loc}
	Let $f\in C^{\infty}_0(\R^2)$ be spherically symmetric, $(x_0,\theta_0)\in T\Sp^1$, 
	$f|_{\gamma{(x_0,\theta_0)}}\not\equiv 0$ and 
	$f|_{\gamma{(x_0,\theta_0)}}$ changes the sign. Then there exist $\varepsilon_0 > 0$ and weight $W_{(x_0, \theta_0), \varepsilon_0}$ such that
	\begin{align}\label{b.loc.weight}
		&P_{W_{(x_0, \theta_0), \varepsilon_0}}f = 0 
		 \text{ on } \Omega(\mathcal{J}_{r(x_0,\theta_0),\varepsilon_0}),  \\
		\begin{split}
		\label{b.loc.weight.props}
		&W_{(x_0, \theta_0),\varepsilon_0} \text{ is bounded, infinitely smooth},
		\\ 
		&\text{strictly positive 
		and rotation invariant on } \Omega (\mathcal{J}_{r(x_0,\theta_0),\varepsilon_0}), 
		\end{split}
	\end{align}
	where $\Omega(\mathcal{J}_{r,\varepsilon_0}), \mathcal{J}_{r,\varepsilon_0}$ are defined 
	in \eqref{loc.weight.domain} and \eqref{int.def}, respectively.
\end{lem}
\par Lemma~\ref{b.lem.loc} is proved in Section~\ref{prv.lm3}. This lemma is a two-dimensional analog of the related lemma in \cite{gonnov2017wrt}.
\begin{rem}
	In Lemma~\ref{b.lem.loc} the construction of $W_{(x_0,\theta_0),\varepsilon_0}$ arises 
	from 
	\begin{enumerate}
		\item finding strictly positive and regular weight $W_{(x_0,\theta_0),\varepsilon}$ on the rays 
	$\gamma = \gamma(x,\theta)$ with fixed $\theta = \theta_0$, where $r(x,\theta_0)\in \mathcal{J}_{r(x_0, \theta_0),\varepsilon}$ for some $\varepsilon > 0$, such that \eqref{b.loc.weight} holds for $\theta = \theta_0$ and 
	under the condition that 
	\begin{equation}
		\hspace{-0.1cm}W_{(x_0,\theta_0),\varepsilon}(y,\gamma) = W_{(x_0,\theta_0),\varepsilon}(|y\theta_0|, \gamma), \, y\in \gamma = \gamma(x,\theta_0), \, r(x,\theta_0)\in \mathcal{J}_{r(x_0,\theta_0),\varepsilon};
	\end{equation}
		\item extending $W_{r(x_0,\theta_0),\varepsilon}$ to all rays $\gamma = \gamma(x,\theta), \, r(x,\theta)\in \mathcal{J}_{r(x_0,\theta_0),\varepsilon}, \, \theta\in \Sp^1$,  via formula \eqref{w_rot_gamma_inv_def}.
	\end{enumerate}
	We recall that $r(x,\theta)$ is defined in \eqref{distance.line.func}.
\end{rem}

\par Let $f$ be the function of \eqref{m.def.ser.}, \eqref{m.def.fk}. Then, using Lemmas~\ref{b.lem.ch.sgn},~\ref{b.lem.loc} one can see that 
\begin{align}\label{exist.loc.weight}
	\begin{split}
	&\text{for all } \delta\in (0,1) \text{ there exist } 
	\{J_i = \mathcal{J}_{r_i,\varepsilon_i}, W_i= W_{(x_i, \theta_i),  \varepsilon_i}\}_{i=1}^N\\
	&\text{such that } J_i, \, i=\overline{1,N}, \text{ is an open cover of } 
	[0,\delta] \text{ in }\R,\\
	&\text{and } W_i \text{ satisfy \eqref{b.loc.weight} and \eqref{b.loc.weight.props} on } \Omega(J_i), \text{ respectively}. 
	\end{split}
\end{align}
\par Actually, we consider  \eqref{exist.loc.weight} for the case of $\delta = \delta_0$ of \eqref{w0.pos.d}.
\par Note that in this case $\{\Omega(J_i)\}_{i=1}^N$ for $J_i$ of \eqref{exist.loc.weight} is the open cover of $\Omega_1(\delta_0)$. 

\par  To the set $\Omega_0(\delta_0)$ we associate the open set
\begin{equation}\label{choice.v0}
J_0= (\delta_0, +\infty)\subset \R.
\end{equation}
Therefore, the collection of intervals $\{\pm J_i, \, 
i = \overline{0,N}\}$ is an open cover of $\R$, where $-J_i$ is the symmetrical reflection of $J_i$ with respect to $\{0\}\in \R$. 
\par We construct the partition of unity $\{\xi_i\}_{i=0}^N$ as follows: 
\begin{align}\label{constr.supp.def}
	&\xi_i(s) = \xi_i(|s|) = \dfrac{1}{2}(\tilde{\xi}_i(s) + \tilde{\xi}_i(-s)), \, 
	s\in \R,\\ \label{constr.supp.prop}
	&\supp \, \xi_i \subset J_i\cup (-J_i), \, i = \overline{0,N},
\end{align}	
where $\{\tilde{\xi}_i\}_{i=0}^{N}$ is a partition of unity for the open cover $\{J_i\cup (-J_i)\}_{i=0}^N$ (see Section~\ref{sect.prelim}, Partition of unity, for $U_i = J_i\cup (-J_i)$).

\par Properties \eqref{supp.xi0}, \eqref{constr.supp.prop} follow from \eqref{gen.part.prop} for $\{\tilde{\xi}_i\}_{i=0}^N$ with $U_i = J_i\cup (-J_i)$, the symmetry of $J_i\cup (-J_i), \,  
i = \overline{1,N}$, choice of $J_0$ in \eqref{choice.v0} and from \eqref{constr.supp.def}.
 
In turn, \eqref{super.prop.xi0} follows from \eqref{choice.v0} and the construction of $J_i, \, 
i = \overline{1,N}$, from \eqref{exist.loc.weight} (see the proof of Lemma~\ref{b.lem.loc} and properties \eqref{exist.loc.weight} in 
Section~\ref{prv.lm3} for details). 
\par Properties \eqref{w.loc.prop}-\eqref{w.loc.prop.last} follow from  \eqref{exist.loc.weight} for $\delta = \delta_0$ and from \eqref{choice.v0}-\eqref{constr.supp.prop}.
\par This completes the description of $W_1, \dots, W_N$ and $\{\xi_i\}_{i=0}^{N}$.

\begin{rem}\label{rem_constr_w_local}
	We have that $J_i = \Lambda_i, \, i = 1,\dots, N$, where $\Lambda_i$ are the intervals 
	in formula \eqref{rem_lambda_mention} of Remark~\ref{rem_general_constr_w} and $J_{i}$ are 
	the intervals considered in \eqref{exist.loc.weight}, \eqref{choice.v0}.
\end{rem}

\subsection{Construction of $W$ and $f$ for $d\geq 3$}\label{subsect.extension}
\par Consider $f$ and $W$ of Theorem~\ref{main.thm}, for $d=2$, constructed in Subsection ~\ref{subsect_construct_dtwo}. For these $f$ and $W$ consider $\tilde{f}$ and $\tilde{U}$
such that 
\begin{equation}\label{two_dim_sym_functions}
	f(x) = \tf(|x|), \, W(x,\theta) = \tilde{U}(|x|, |x\theta|), \, x\in \R^2, \, \theta
	\in \Sp^1.
\end{equation}
\begin{prop}\label{corollary.generalizations}
	Let $W$ and $f$, for $d\geq 3$, be defined as
	\begin{align}\label{generalization.w}
		&W(x,\theta) = \tilde{U}(|x|, |x\theta|), \, (x,\theta)\in \R^d\times \Sp^{d-1},\\
		\label{generalization.f}
		&f(x) = \tilde{f}(|x|), \, x\in \R^d,
	\end{align}
	where $\tilde{U}, \, \tilde{f}$ are the functions of \eqref{two_dim_sym_functions}.
	Then
	\begin{equation}\label{high_dim_extension}
		P_Wf \equiv 0 \text{ on } T\Sp^{d-1}.
	\end{equation}
	In addition, weight $W$ satisfies properties \eqref{thm_prop_smooth}-\eqref{thm_w_unity}, 
	$f$ is spherically symmetric infinitely smooth and compactly supported on $\R^d$, $f\not\equiv 0$.
\end{prop}
\par Proposition~\ref{corollary.generalizations} is proved in Subsection~\ref{proof.corollary.generalizations}.
\par This completes the proof of Theorem~\ref{main.thm}.

\section{Proof of Theorem~\ref{corol_dimension_infty}}\label{sect_thm2_proof}
\subsection{Proof for $d\geq 3$}
\par Let 
\begin{align}
	\label{thm2_proof_Wdef}
	&W \text{ be the weight of Theorem~\ref{main.thm} for }
	d\geq 3, \\ \label{thm2_proof_Rdef}
	&R \text{ be the number in \eqref{thm_w_unity} for } d\geq 3, \\
	\label{thm2_proof_yidef}
	\begin{split}
	&\{y_i\}_{i=1}^{\infty} \text{ be a sequence of vectors in } \R^d 
	\text{ such that } y_1 = 0, \, |y_i - y_j| > 2R \\
	&\text{for } i\neq j,\, i,j\in \mathbb{N}, 
	\end{split}\\ 
	\label{thm2_proof_bidef}
	&\{\overline{B}_i\}_{i=1}^{\infty} \text{ be the closed balls in } \R^d \text{ of radius } R 
	\text{ centered at } y_i\,  (\text{see }\eqref{thm2_proof_Rdef}, \eqref{thm2_proof_yidef}).
\end{align}
The weight $W_n$ is defined as follows
\begin{align}\label{thm2_proof_wndef}
	&W_n(x,\theta) = \begin{cases}
		&1 \text{ if } x\not\in \bigcup\limits_{i=1}^n \overline{B}_i, \\
		&W(x-y_1,\theta) = W(x, \theta) \text{ if } x\in \overline{B}_1, \\
		&W(x-y_2, \theta) \text{ if } x\in \overline{B}_2,\\
		&...,\\
		&W(x-y_k, \theta) \text{ if } x\in \overline{B}_k,\\
		&...,\\
		&W(x-y_n, \theta) \text{ if } x\in \overline{B}_n,
	\end{cases}\\ \nonumber
	&\theta\in \Sp^{d-1}, \, n\in \mathbb{N}\cup \{\infty\}, \, d\geq 3,
\end{align}
where $W$ is defined in \eqref{thm2_proof_Wdef}, $y_i$ and $\overline{B}_i$ are defined in \eqref{thm2_proof_yidef}, \eqref{thm2_proof_bidef}, respectively.
\par Properties \eqref{w.compl.},  \eqref{wn_smooth_d3} and \eqref{wn_identity_outside} for $W_n$, defined in \eqref{thm2_proof_wndef}, for $d\geq 3$, follow from \eqref{thm_prop_smooth}-\eqref{thm_w_unity}, \eqref{thm2_proof_Wdef}, \eqref{thm2_proof_Rdef}.
\par Let
\begin{align}
\begin{split}\label{thm2_proof_fidef}
	&f_1(x) \eqdef f(x), \, f_2(x) \eqdef f(x-y_2),\dots, f_n(x) \eqdef f(x-y_n), \, x\in \R^d, \, d\geq 3,
\end{split}
\end{align}
where $y_i$ are defined in \eqref{thm2_proof_yidef} and 
\begin{equation}\label{thm2_proof_fidefthm}
	f \text{ is the function of Theorem~\ref{main.thm} for } d\geq 3.
\end{equation}
One can see that 
\begin{equation}\label{thm2_proof_fiprops}
	f_i\in C_0^{\infty}(\R^d),\, d\geq 3, \,  f_i \not \equiv 0, \,  \supp \, f_i \subset \overline{B}_i, \, 
	\overline{B}_i\cap \overline{B}_j = \emptyset \text{ for } i\neq j,
\end{equation}
where $\overline{B}_i$ are defined in \eqref{thm2_proof_bidef}, 
$i = 1, \dots, n$.
\par The point is that 
\begin{align}\label{thm2_proof_fiinkernel}
	&P_{W_n}f_i\equiv 0 \text{ on }T\Sp^{d-1}, \, 
	d\geq 3, \, i = 1, \dots, n, \\ 
	\label{thm2_proof_fiindep}
	&f_{i} \text{ are linearly independent in } C_0^{\infty}(\R^d), \, d\geq 3, 
	\, i = 1, \dots, n,
\end{align}
where $W_n$ is defined in \eqref{thm2_proof_wndef}, $f_i$ are defined in \eqref{thm2_proof_fidef}. 

To prove \eqref{thm2_proof_fiinkernel} we use, in particular, the following general formula:
\begin{align}\label{thm2_proof_pwgenprop}
	\begin{split}
		&P_{W_y}f_y(x,\theta) = 
		\int\limits_{\gamma(x,\theta)}W(y'-y,\theta)f(y'-y)dy' \\
		&\qquad \qquad \quad = 
		\int\limits_{\gamma(x-y, \theta)}W(y',\theta)f(y')dy' 
		= P_Wf(x-y, \theta), \, x\in \R^d, \, \theta\in \Sp^{d-1},
	\end{split}\\
	\label{thm2_proof_wfshift}
	&W_y(x,\theta) = W(x-y,\theta), \, f_y = f(x-y), \, 
		 \, x,\, y\in \R^d, \, \theta\in \Sp^{d-1}.
\end{align}
where $W$ is an arbitrary weight satisfying \eqref{w.compl.}, $f$ is a test-function, 
$\gamma(x,\theta)$ is defined according to \eqref{ray_def_omega}.
\par Formula \eqref{thm2_proof_fiinkernel} follows from formula \eqref{m.thm.stm.}, definitions \eqref{thm2_proof_wndef}, \eqref{thm2_proof_fidef}, \eqref{thm2_proof_fidefthm}, 
properties \eqref{thm2_proof_fiprops} and from formulas \eqref{thm2_proof_pwgenprop}, \eqref{thm2_proof_wfshift}.
\par Formula \eqref{thm2_proof_fiindep} follows from definitions \eqref{thm2_proof_fidef}, \eqref{thm2_proof_fidefthm} and properties \eqref{thm2_proof_fiprops}. 
\par This completes the proof of Theorem~\ref{corol_dimension_infty} for $d\geq 3$.

\subsection{Proof for $d= 2$}
In \cite{boman1993example}, there were constructed a weight $W$ and a function $f$ for $d=2$, such that:
\begin{align}\label{thm2_proof_pwbzero}
&P_{W}f \equiv 0 \text{ on }T\Sp^1, \\
\label{thm2_proof_wbpos}
&W = \overline{W} \geq c > 0, \, W\in C^{\infty}(\R^2\times \Sp^{1}),\\
\label{thm2_proof_fbprops}
&f\in C_0^{\infty}(\R^2), \, f \not \equiv 0, \, \supp f\subset \overline{B},	
\end{align}
where $c$ is a constant, $\overline{B}$ is defined in \eqref{ball_def}.
\par We define 
\begin{equation}\label{thm2_proof_wd2def}
	\tW(x,\theta) = c^{-1}\phi_1(x)W(x,\theta) + \phi_2(x), \, x\in \R^2, \, \theta\in \Sp^{1},
\end{equation}
where $W$ is the weight of \eqref{thm2_proof_pwbzero}, \eqref{thm2_proof_wbpos}, $c$ is a 
constant of \eqref{thm2_proof_wbpos}.
\begin{align}\label{thm2_proof_unity_def}
	\begin{split}
&\phi_1 = \phi_1(x), \phi_2 = \phi_2(x) \text{ is a $C^{\infty}$-smooth partition of unity on } \R^2 \text{ such that}, \\
&\phi_1 \equiv 0 \text{ for } |x|\geq R > 1, \, \phi_1 \equiv 1 \text{ for } |x|\leq 1, 
\, \phi_1 \geq 0 \text{ on } \R^2,\\
&\phi_2 \equiv 0 \text{ for } |x|\leq 1, \, \phi_2\geq 0 \text{ on }\R^2, 
\end{split}
\end{align}
where $R$ is a constant.
\par From \eqref{thm2_proof_pwbzero}-\eqref{thm2_proof_unity_def}  it follows that 
\begin{align}\label{thm2_proof_twprops1}
	&P_{\tW}f \equiv 0 \text{ on } T\Sp^1, \\
	\label{thm2_proof_twprops2}
	\begin{split}
	&\tW \geq 1, \,  \tW \in C^{\infty}(\R^2\times \Sp^1),\\
	&\tW(x,\theta) \equiv 1 \text{ for } |x| \geq R > 1, \, \theta\in \Sp^1.
	\end{split}
\end{align}
\par The proof of Theorem~\ref{corol_dimension_infty} for $d= 2$ proceeding from 
\eqref{thm2_proof_fbprops}, 
\eqref{thm2_proof_wd2def}, \eqref{thm2_proof_twprops1}, \eqref{thm2_proof_twprops2} is completely 
similar to the proof of Theorem~\ref{corol_dimension_infty} for $d\geq 3$, proceeding 
from Theorem~\ref{main.thm}.
\par Theorem~\ref{corol_dimension_infty} is proved.

\section{Proofs of Corollary~\ref{thm_delta_non_inj} and Proposition~\ref{corollary.generalizations}}\label{proof.coroll.prop1}
\subsection{Proof of Corollary~\ref{thm_delta_non_inj}}\label{proof_cor1}
	Let 
	\begin{align}
		&X_r = \{x_1e_1 + x_2e_2 + r e_3 : (x_1,x_2) \in \R^2\}, \,   0 \leq r < 1,\\
		&S = X_0\cap \Sp^{2} = \{(\cos\phi, \sin\phi, 0) \in \R^{3} : 
		\phi \in [0,2\pi)\} \simeq \Sp^1.
	\end{align}	 
	where $(e_1,e_2,e_3)$ is the standard orthonormal basis in $\R^3$.
	\par Without loss of generality we assume that $0 < \delta < 1$.
	Choosing $r$ so that \newline ${\sqrt{1-\delta^2} \leq  r < 1}$, we have that the intersection of the three dimensional ball $B(0,1)$ with $X_r$ is the two-dimensional disk $B(0,\delta'), \, 
	\delta' \leq \delta$ (with respect to the coordinates $(x_1,x_2)$ induced by basis $(e_1,e_2)$ on $X_r$). 
	\par We define $N, \, W_\delta$ on $\R^2\times \Sp^1$ and $f_\delta$ on $\R^2$ 
	as follows:
	\begin{align}
		\label{N_norm_def}
		&N = \|W\|_{C^{\alpha}(\R^3\times\Sp^2)},\\
		\label{weight_restrictions}
		&W_\delta := W|_{X_{r}\times S},\\
		\label{func_restriction}
		&f_\delta := f|_{X_{r}},\\ \nonumber
		&\text{for } r=\sqrt{1-\delta^2},
	\end{align}
	where $W$ and $f$ are the functions of Theorem~\ref{main.thm} for $d=3$.
	\par Due to \eqref{thm_prop_smooth}-\eqref{thm_prop_pos}, 
	\eqref{N_norm_def},  \eqref{weight_restrictions} we have that 
	\begin{equation}\label{wdelta_props}
		W_\delta \geq 1/2, \, \|W_\delta\|_{C^{\alpha}(\R^2\times \Sp^1)} \leq N.
	\end{equation}
	\par Properties \eqref{wdelta_props} imply \eqref{cexample_inj2_start}. 
	\par In view of Lemma~\ref{b.lem.ch.sgn} for the function $f$ of Theorem~\ref{main.thm}, we have that $f_\delta$ is spherically symmetric, $f_\delta\in C^{\infty}_0(B(0,\delta')), \, 
	f_\delta\not\equiv 0$.	
	\par Using \eqref{m.thm.stm.}, \eqref{weight_restrictions}, \eqref{func_restriction} one can see that \eqref{cexample_inj2_end} holds.
	\par This completes the proof of Corollary~\ref{thm_delta_non_inj}.

\subsection{Proof of Proposition~\ref{corollary.generalizations}}\label{proof.corollary.generalizations}
\par Let 
\begin{equation}\label{extension_Ir_def}
	I(r) = 
	\int\limits_{\gamma_r}
	\tilde{U}(|y|, r)\tilde{f}(|y|) \, dy, \, r\geq 0, \, 
	\gamma_r = \gamma(re_2,e_1),
\end{equation}
where $\gamma(x,\theta)$ is defined by \eqref{ray_in_Rd}, $(e_1,\dots ,e_d)$ is the standard 
basis in $\R^d$.
\par Due to formula \eqref{m.thm.stm.} of Theorem~\ref{main.thm} for $d=2$ and 
formulas \eqref{two_dim_sym_functions}, \eqref{extension_Ir_def} we have that 
\begin{equation}\label{extension_Ir_zero}
I(r) = P_Wf(re_2, e_1) = 0 \text{ for } r\geq 0.
\end{equation}
Next, using \eqref{rad.def}, \eqref{two_dim_sym_functions}, \eqref{extension_Ir_zero} we have also that 
\begin{align}\label{extension_zero_ident}
	&P_Wf(x,\theta) = \int\limits_{\gamma(x,\theta)}
	\tilde{U}(|y|, |y-(y\theta)\theta|)\tilde{f}(|y|)\, dy = I(|x|) = 0 
	\text{ for } (x,\theta)\in T\Sp^{d-1},
\end{align}
where $\gamma(x,\theta)$ is defined in \eqref{ray_in_Rd}.
\par Formula \eqref{extension_zero_ident} implies \eqref{high_dim_extension}.  Properties of $W$ and $f$ mentioned in Proposition~\ref{corollary.generalizations} follow from 
properties \eqref{thm_prop_smooth}-\eqref{thm_w_unity} of $W$ and of $f$ of Theorem~\ref{main.thm} for $d=2$.
\par This completes the proof of Proposition~\ref{corollary.generalizations}.

\section{Proofs of formulas \eqref{w0.c.rinv}-\eqref{w0.zero}}
\label{proofs.manyformulas.w0}

\subsection{Proof of formulas \eqref{w0.c.rinv}-\eqref{w0.pos.d}}\label{subsect_pr_lem2}
\begin{lem}\label{b.cont.lem}
	Let $W_0$ be defined by \eqref{b.w0.def}, \eqref{b.gkhk.def}. Then $W_0$ admits the following representation:
		\begin{align}\label{w_0_eq_U0}
			&W_0(x,\theta) = U_0(x\theta, |x-(x\theta)\theta|), \, 
			(x,\theta)\in \Omega((1/2, +\infty)),\\
			\label{U0_def}
			&U_0(s,r) = \begin{cases}
				&1 - \widetilde{G}(r)\sum\limits_{k=3}^{\infty}k! \tf_k((s^2 + r^2)^{1/2})
				\dfrac{\psi_{k-2}(r)}{\widetilde{H}_k(r)}, \, 1/2 < r < 1, \\
				&1, \, r \geq 1
			\end{cases}, \\ \label{ghk_def}
			\begin{split}
			&\tG(r) \eqdef \int\limits_{\gamma_r}
			\tf(|y|)\,dy, \, \tH_k(r) \eqdef \int\limits_{\gamma_r}\tf_k^2(|y|)\, dy,
			\, \tf = \sum\limits_{k=1}^{\infty}\dfrac{\tf_k}{k!},
			\end{split}
			\\ \nonumber
			&s\in \R, \, x\in \R^2, \, \gamma_r \text{ is an arbitrary ray in } 
			T(r), \, r > 1/2,
		\end{align}
		where $\tf_k$ are defined by \eqref{m.def.fk}, $T(r)$ is defined by \eqref{t_def}, $d=2$. 
		In addition:
			\begin{align}\label{U0_smooth_property}
				&U_0 \text{ is infinitely smooth on } \R \times \{(1/2,1)\cup (1,+\infty)\},\\
			\label{u0.border.convergence}
			&U_0(s,r) \rightarrow 1 \text{ as }r\rightarrow 1 
			\text{ (uniformly in } s\in \R\text{)},\\ 
			\label{u0.identity.exterior}
			&U_0(s,r) = 1 \text{ if } s^2 + r^2 \geq 1,\\ 		
			\label{w.w0.conv}
			&|1-U_0(s,r)| \leq C_0(1-r)^{1/2}
				\log^4_2 \left( \dfrac{1}{1-r}\right),\\ \nonumber
			&\hspace{1.0cm} \text{for } s\in \R, \, 1/2 < r < 1,\\
			\label{U0_holder_estimate}
			&|U_0(s,r) - U_0(s',r')| \leq C_1 
	 			|s-s'|^{1/\alpha} + C_1 |r - r'|^{1/\alpha},\\ \nonumber
	 		&\hspace{1.0cm} \text{for } \alpha\in (0,1/16), \, s, s'\in \R, \, 
	 		r,r' > 1/2,
		\end{align}
			where $C_0, C_1$ are positive constants depending on $\Phi$ of
	 \eqref{phi.suppint}-\eqref{phi.unity.cond}.
\end{lem}
\par Lemma~\ref{b.cont.lem} is proved Section~\ref{proof_bcontlem}.
\par Lemma~\ref{b.cont.lem} implies \eqref{w0.c.rinv}-\eqref{w0.pos.d} as follows.
\par The continuity and rotation invariancy of $W_0$ in \eqref{w0.c.rinv} follow from \eqref{w.rot.inv}, \eqref{u.rot.inv.prop}, \eqref{w_0_eq_U0}, \eqref{U0_holder_estimate}. 
\par Due to \eqref{psik.unit.part.def}, \eqref{w_0_eq_U0}, \eqref{U0_def}, 
\eqref{ghk_def} we have also that 
\begin{equation}\label{U0_left_extension}
U_0 \text{ admits a continuous extension to } \R\times [1/2, +\infty).
\end{equation}
Properties \eqref{u0.identity.exterior}, \eqref{U0_left_extension} imply the boundedness of $W_0$ on $\Omega_0(1/2)$, where $\Omega_0(\cdot)$ is defined in \eqref{pr.l0}, $d=2$. This completes the proof of \eqref{w0.c.rinv}.
\par Formula \eqref{w0_holder_property} follows from \eqref{w_0_eq_U0}, \eqref{U0_smooth_property}, \eqref{U0_holder_estimate} and from 
the fact that $x\theta$, ${|x-(x\theta)\theta|}$ are infinitely smooth functions on $\Omega_0(1/2)$ 
and are Lipshitz in $(x,\theta)$ for ${x\in \overline{B}(0,R), \, R>1}$.
\par Formula \eqref{w0.pos.d} follows from \eqref{w0.c.rinv}, \eqref{w_0_eq_U0}, \eqref{U0_def}, 
 \eqref{u0.border.convergence}, \eqref{u0.identity.exterior}.
\par This completes the proof of \eqref{w0.c.rinv}-\eqref{w0.pos.d}.

\subsection{Proof of formula \eqref{w0.zero}}\label{subs.w0.zero}
\par From \eqref{rad.def}, \eqref{m.def.ser.}-\eqref{phi.supp}, \eqref{psik.unit.part.def}, \eqref{b.w0.def}, \eqref{b.gkhk.def} it follows that:
\begin{align}\nonumber
	P_{W_0}f(x,\theta) &= \int\limits_{\gamma(x,\theta)}
	\hspace*{-0.2cm}f(y)\, dy - 
	G(x,\theta)
	\sum\limits_{k=3}^{\infty}k!\psi_{k-2}(r(x,\theta))
	\dfrac{\int\limits_{\gamma(x,\theta)}f(y)f_k(y)dy}
	{H_k(x,\theta)} \\
	&= \int\limits_{\gamma(x,\theta)}\hspace*{-0.2cm}f(y)\, dy - 
	\int\limits_{\gamma(x,\theta)}\hspace*{-0.2cm}f(y)\, dy
	\sum\limits_{k=3}^{\infty}\psi_{k-2}(r(x,\theta))
	\dfrac{\int\limits_{\gamma(x,\theta)}f_k^2(y)dy}
	{\int\limits_{\gamma(x,\theta)}f_k^2(y)\,dy}\\ \nonumber
	&= \int\limits_{\gamma(x,\theta)}\hspace*{-0.2cm}f(y)\, dy - 
	\int\limits_{\gamma(x,\theta)}\hspace*{-0.2cm}f(y)\, dy
	\sum\limits_{k=3}^{\infty}\psi_{k-2}(r(x,\theta)) = 0 \text{ for } (x,\theta)\in \Omega_0(1/2),
\end{align}
where $\gamma(x,\theta)$ is defined in \eqref{ray_in_Rd}, $\Omega_0(\cdot)$ is defined in \eqref{pr.l0}, $d=2$.
\par Formula \eqref{w0.zero} is proved.

\section{Proof of Lemma~\ref{b.lem.loc}}\label{prv.lm3}
\par By $u\in \R$ we denote the coordinates on a fixed ray $\gamma(x,\theta), \, 
(x,\theta)\in \Omega, \, d=2$, taking into account the orientation, where $u = 0$ at the point $x-(x\theta)\theta\in \gamma(x,\theta)$; see notation \eqref{ray_def_omega}. 

\par Using Lemma~\ref{b.lem.ch.sgn}, one can see that
\begin{equation}\label{plane.symm}
	f|_{\gamma{(x,\theta)}}\in C_0^{\infty}(\R), \, f|_{\gamma{(x,\theta)}}(u) = f|_{\gamma{(x,\theta)}}(|u|), \, u\in \R.
\end{equation}
\par Using \eqref{plane.symm} and the assumption that  
$f|_{\gamma{(x_0,\theta_0)}}(u)$ changes the sign, one can see that there exists $\psi_{(x_0,\theta_0)}$ such that
\begin{align}\label{psi.loc.symm}
	&\psi_{(x_0,\theta_0)}\in C_0^{\infty}(\R), 
	\, \psi_{(x_0,\theta_0)} \geq 0, \,  
	\psi_{(x_0,\theta_0)}(u) = \psi_{(x_0,\theta_0)}(|u|), \,  u\in \R,\\
	\label{psi_loc_int_non_zero}
	&\int\limits_{\gamma(x_0,\theta_0)} \hspace{-0.3cm}
	f\psi_{(x_0,\theta_0)} \, d\sigma \neq 0,
\end{align}
and if 
\begin{equation}
	\int\limits_{\gamma{(x_0,\theta_0)}}\hspace*{-0.3cm}f\,d\sigma \neq 0
\end{equation}
then also 
\begin{align}
	\label{sgn.choice}
	&\sgn
	(
	\hspace*{-0.3cm}
		\int\limits_{\gamma{(x_0,\theta_0)}} \hspace*{-0.35cm} f\, d\sigma
	)\, 
	\sgn
	(	\hspace*{-0.3cm}
		\int\limits_{\gamma{(x_0,\theta_0)}} \hspace*{-0.35cm} f\psi_{(x_0,\theta_0)}\, d\sigma 
	) = -1,
\end{align}
where $d\sigma = du$ (i.e., $\sigma$ is the standard Euclidean measure on $\gamma{(x,\theta)}$).
\par Let 
\begin{equation}\label{loc.weight}
	W_{(x_0, \theta_0)}(x,\theta) = 1-\psi_{(x_0,\theta_0)}(x\theta)
	\dfrac{\int\limits_{\gamma{(x,\theta)}}\hspace{-0.3cm}
	f\, d\sigma}{\int_{\gamma{(x,\theta)}}\hspace{-0.1cm}f\psi_{(x_0,\theta_0)}\, d\sigma},
	\, x\in \R^2, \, \theta\in \Sp^1,
\end{equation}
where $d\sigma = du$, where $u$ is the coordinate on $\gamma(x,\theta)$.
\par Lemma~\ref{b.lem.ch.sgn} and property \eqref{psi.loc.symm} imply that
\begin{equation}\label{save.sign}
	\int\limits_{\gamma{(x,\theta)}}\hspace{-0.3cm}
	f\, d\sigma \text{ and }\hspace{-0.2cm}
	\int\limits_{\gamma{(x,\theta)}}\hspace{-0.3cm}
	f\psi_{(x_0,\theta_0)}\, d\sigma \text{ depend only on }
	r(x,\theta), \text{ where } (x,\theta)\in \Omega, \,
\end{equation}
where $r(x,\theta)$ is defined in \eqref{distance.line.func}, $\Omega$ is defined in \eqref{pr.ts2}, $d=2$. 
\par From \eqref{psi.loc.symm}, \eqref{loc.weight}, \eqref{save.sign} it follows that $W_{(x_0, \theta_0)}$ is rotation-invariant in the sense \eqref{w.rot.inv}.
\par Formulas \eqref{psi_loc_int_non_zero}, \eqref{loc.weight}, \eqref{save.sign}, properties of $f$ of Lemma~\ref{b.lem.ch.sgn} and properties of $\psi_{(x_0,\theta_0)}$ of  \eqref{psi.loc.symm} imply that 
\begin{align}\label{ws0.contin0}
	\exists \varepsilon_1 > 0 : \int\limits_{\gamma{(x,\theta)}}\hspace*{-0.3cm}f\psi_{(x_0,\theta_0)}\,d\sigma \neq 0 \text{ for }  
	(x,\theta)\in \Omega(\mathcal{J}_{r(x_0,\theta_0),\varepsilon_1}),
\end{align} 
where sets $\Omega(\mathcal{J}_{s,\varepsilon}), \, \mathcal{J}_{s,\varepsilon}$ are defined in 
\eqref{loc.weight.domain}, \eqref{int.def}, respectively.
\par In addition, using properties of $f$ of Lemma~\ref{b.lem.ch.sgn} and also using 
\eqref{m.def.ser.}, \eqref{f_k_supports}, \eqref{psi.loc.symm}, \eqref{loc.weight}, \eqref{ws0.contin0}, one can see that
\begin{align}\label{ws0.contin}
	W_{(x_0, \theta_0)}\in C^{\infty}(\Omega(\mathcal{J}_{r(x_0,\theta_0),\varepsilon_1})).
\end{align}
\par In addition, from \eqref{plane.symm}-\eqref{save.sign} it follows that 
\begin{align}\nonumber
	\text{ if } r(x,\theta) = r(x_0,\theta_0) \text{ then } W_{(x_0, \theta_0)}(x,\theta) &= 1-\psi_{(x_0,\theta_0)}(x\theta)
	\dfrac{\int\limits_{\gamma{(x_0,\theta_0)}}\hspace{-0.4cm}
	f\, d\sigma}
	{\int\limits_{\gamma{(x_0,\theta_0)}}\hspace{-0.4cm}
	f\psi_{(x_0,\theta_0)}\, d\sigma} \\ \label{origin.pos}
	&= 1-\psi_{(x_0,\theta_0)}(x\theta)
	\dfrac{\int\limits_{\gamma{(x_0,\theta_0)}}\hspace{-0.4cm}
	f\, d\sigma}{\int\limits_{\gamma{(x_0,\theta_0)}}\hspace{-0.4cm}
	f\psi_{(x_0,\theta_0)}\, d\sigma} \geq 1,
\end{align}
where $r(x,\theta)$ is defined in \eqref{distance.line.func}, $d=2$.

From properties of $f$ of Lemma~\ref{b.lem.ch.sgn}, properties of $\psi_{(x_0,\theta_0)}$ of \eqref{psi.loc.symm} and from formulas \eqref{loc.weight}, \eqref{ws0.contin0}, \eqref{ws0.contin}, \eqref{origin.pos} it follows that 

\begin{equation}\label{pos.cond.zero}
	\exists \varepsilon_0 > 0\, (\varepsilon_0 < \varepsilon_1) : W_{(x_0,\theta_0)}(x,\theta) \geq 1/2 \text{ for }
	(x,\theta)\in \Omega(\mathcal{J}_{r(x_0,\theta_0),\varepsilon_0}).
\end{equation}
Let 
\begin{equation}\label{loc_weight_zero_def}
	W_{(x_0,\theta_0),\varepsilon_0}: = W_{(x_0,\theta_0)} \text{ for } 
	(x,\theta)\in \Omega(\mathcal{J}_{r(x_0,\theta_0),\varepsilon_0}),
\end{equation}
where $W_{(x_0,\theta_0)}$ is defined in \eqref{loc.weight}.
\par Properties \eqref{save.sign}, \eqref{ws0.contin}, \eqref{pos.cond.zero} imply 
\eqref{b.loc.weight.props} for $W_{(x_0, \theta_0),\varepsilon_0}$ of \eqref{loc_weight_zero_def}.
\par Using \eqref{rad.def}, \eqref{loc.weight}, \eqref{ws0.contin0}, \eqref{loc_weight_zero_def}  one can see that
\begin{align} \nonumber
	P_{W_{(x_0,\theta_0), \varepsilon_0}}f(x,\theta) &= \int\limits_{\gamma{(x,\theta)}}\hspace{-0.2cm}
	W_{(x_0,\theta_0)}(\cdot, \theta)f\, d\sigma
	\\ \label{loc_weight_zero_ident}
	&=\int\limits_{\gamma{(x,\theta)}}\hspace{-0.2cm}f\, d\sigma \, - 
	\dfrac{\int\limits_{\gamma{(x,\theta)}}\hspace{-0.3cm}
	f\, d\sigma}{\int\limits_{\gamma{(x,\theta)}}\hspace{-0.3cm}
	f\psi_{(x_0,\theta_0)}\, d\sigma}	
	\int\limits_{\gamma{(x,\theta)}}\hspace{-0.3cm} f \psi_{(x_0,\theta_0)}\, d\sigma = 0 
	\, \text{ for } (x,\theta)\in \Omega(\mathcal{J}_{r(x_0,\theta_0),\varepsilon_0}),
\end{align}
where $\Omega(\cdot)$ is defined in \eqref{loc.weight.domain}, $d=2$, $\mathcal{J}_{r,\varepsilon}$ is defined in \eqref{int.def}.
Formula \eqref{b.loc.weight} follows from \eqref{loc_weight_zero_ident}.
\par Lemma~\ref{b.lem.loc} is proved.

\section{Proof of Lemma~\ref{b.cont.lem}}\label{proof_bcontlem}

\paragraph{\large Proof of \eqref{w_0_eq_U0}-\eqref{ghk_def}.} Using \eqref{distance.line.func}, \eqref{m.def.ser.}, 
\eqref{m.def.fk}, \eqref{b.gkhk.def}, \eqref{ghk_def} we obtain
\begin{align}\label{G_left_right}
	&G(x,\theta) = \tG(r(x,\theta)) = \hspace{-0.3cm}
	\int\limits_{\gamma(x,\theta)}
	\hspace{-0.3cm}f(x)\, dx, \\
	\label{Hk_left_right}
	&H_k(x,\theta) = \tH_k(r(x,\theta)) =  \hspace{-0.3cm}
	\int\limits_{\gamma(x,\theta)}
	\hspace{-0.3cm}f_k^2(x) \, dx,\\
	\label{f_k_left_right}
	&\tf_k(|x|) = \tf_k((|x\theta|^2 + |x-(x\theta)\theta|^2)^{1/2}), \, 
	(x,\theta)\in \Omega_0(1/2),
\end{align}
where $\Omega_0(\cdot)$ is defined in \eqref{pr.l0}, $d=2$, $\gamma(x,\theta)$ is defined as in \eqref{ray_def_omega}.
\par Formulas \eqref{b.w0.def}, \eqref{b.gkhk.def},  \eqref{G_left_right}-\eqref{f_k_left_right} imply \eqref{w_0_eq_U0}-\eqref{ghk_def}.

\paragraph{\large Proof of \eqref{U0_smooth_property}.}
Let 
\begin{equation}\label{lambdak_def}
	\Lambda_k = (1-2^{-k+3}, 1-2^{-k+1}), \, k\in \mathbb{N}, \, k\geq 4.
\end{equation}
From \eqref{psik.unit.part.def} it follows that, for $k\geq 4$:
\begin{align}\label{item2_supp.psi.1}
	&\supp \, \psi_{k-1} \subset  (1-2^{-k+2}, 1-2^{-k}),\\ 
	&\supp \, \psi_{k-2} \subset  (1-2^{-k+3}, 1-2^{-k+1}) = \Lambda_k, \\ 
	\label{item2_supp.psi.3}
	&\supp \, \psi_{k-3} \subset  (1-2^{-k+4}, 1-2^{-k+2}). 
\end{align}
Due to \eqref{U0_def}, \eqref{ghk_def},
 \eqref{item2_supp.psi.1}-\eqref{item2_supp.psi.3}, we have the following formula for $U_0$:
\begin{align}\label{item2_U0_long}\nonumber
		U_0(s,r) = 1 - \tG(r)
		&\left(
			(k-1)!\tf_{k-1}((s^2+r^2)^{1/2})\dfrac{\psi_{k-3}(r)}{\tH_{k-1}(r)}  \right. \\ 
		&\left. 
		 	+k!\tf_{k}((s^2+r^2)^{1/2})\dfrac{\psi_{k-2}(r)}{\tH_{k}(r)} \right. \\
		 	\nonumber
		&\left. 
			+(k+1)!\tf_{k+1}((s^2+r^2)^{1/2})\dfrac{\psi_{k-1}(r)}{\tH_{k+1}(r)}
		\right) \text{ for } r\in \Lambda_k, \, s\in \R, \, k\geq 4.
\end{align}
From \eqref{ghk_def}, \eqref{item2_U0_long} it follows that
\begin{alignat}{3}\nonumber
	&\dfrac{\partial^{n}U_0}{\partial s^n}(s,r) = &&-\tG(r)
	\left(
			(k-1)!
			\dfrac{\partial^n \tf_{k-1}((s^2+r^2)^{1/2})}{\partial s^n}
			\dfrac{\psi_{k-3}(r)}{\tH_{k-1}(r)}  \right. \\  \label{U0_p_high_deriv}
		 &  &&\left. 
		 	+k!
		 	\dfrac{\partial^n \tf_{k}((s^2+r^2)^{1/2})}{\partial s^n}
		 	\dfrac{\psi_{k-2}(r)}{\tH_{k}(r)} \right. \\ \nonumber
		 & &&\left. 
			+(k+1)!
			\dfrac{\partial^n \tf_{k+1}((s^2+r^2)^{1/2})}{\partial s^n}
			\dfrac{\psi_{k-1}(r)}{\tH_{k+1}(r)}
		\right),\\ \label{U0_r_high_deriv}
	&\dfrac{\partial^n \tG}{\partial r^n}(r) =\int\limits_{-\infty}^{+\infty} &&
	\dfrac{\partial^n }{\partial r^n}\tf((s^2 + r^2)^{1/2})\, ds, \, 
	\dfrac{\partial^n \tH_m}{\partial r^n}(r) = \int\limits_{-\infty}^{+\infty}
	\dfrac{\partial^n}{\partial r^n}
	\tf_m^2((s^2 + r^2)^{1/2})\, ds, 	\\ \nonumber
	& &&
	\hspace*{-1.9cm}
	 \, 
	r\in \Lambda_k, \, s \in  \R,
	\, m \geq 1, 
	\, n \geq 0, 
	\, k\geq 4,
\end{alignat}
where $\tG, \, \tH_m$ are defined in \eqref{ghk_def}.
\par Using Lemma~\ref{b.lem.ch.sgn} and formulas \eqref{m.def.ser.}, \eqref{m.def.fk}, \eqref{psik.unit.part.def}-\eqref{H_k_smoothness}, \eqref{ghk_def} one can see that: 
\begin{align}\label{item2_smoothness_fact}
\begin{split}
	&\tf, \, \tf_{m-2}, \, \tG, \, \tH_m \text{ belong to } C_0^{\infty}(\R),\\ 
 	&\dfrac{\psi_{m-2}}{\tH_m} \text{ belongs to }
	C^{\infty}_0((1/2,1)) \text{ for any } m\geq 3.
\end{split}
\end{align}
\par From \eqref{U0_p_high_deriv}-\eqref{item2_smoothness_fact} it follows that
$U_0(s,r)$ has continuous partial derivatives of all orders with respect to $r\in \Lambda_k, \, s\in \R$. It implies that $U_0\in C^{\infty}(\R\times \Lambda_k)$.
From the fact that $\Lambda_k, \, k\geq 4$, is an open cover of $(1/2, 1)$ and 
from definition \eqref{U0_def} of $U_0$, it follows that $U_0\in C^{\infty}(\R\times \{(1/2,1)\cup (1,+\infty)\})$. 
\par This completes the proof of \eqref{U0_smooth_property}.

\paragraph{ \large Proof of \eqref{u0.identity.exterior}.}
\par From \eqref{m.def.fk}-\eqref{phi.unity.cond} it follows that  
\begin{equation}\label{fk_support_bound}
	\tf_k(|x|) = 0 \text{ if } |x|\geq 1 \text{ for } k\in \mathbb{N}. 
\end{equation}
Formula $|x|^2 = |x\theta|^2 + |x-(x\theta)\theta|^2, \, x\in \R^2, \, \theta\in \Sp^1$, 
and formulas \eqref{U0_def}, \eqref{fk_support_bound} imply \eqref{u0.identity.exterior}.

\paragraph{\large Proofs of \eqref{w.w0.conv}-\eqref{U0_holder_estimate}.}
\begin{lem}\label{main_lemm}
There are positive constants $c, k_1$ depending on $\Phi$ of \eqref{phi.suppint}-\eqref{phi.unity.cond}, such that
\begin{itemize}
	\item[(i)] for all $k\in \mathbb{N}$ the following estimates hold:
	\begin{align}\label{estimate_func}
		&|\tf_k| \leq 1,\\
		\label{estim_Dfunc}
		&|\tf_k'| \leq c 8^k,
	\end{align}
	where $\tf_k'$ denotes the derivative of $\tf_k$ defined in \eqref{ghk_def}.
	\item[(ii)] for $k\geq k_1$ and $1/2 < r \leq 1$ the following estimates hold:
	\begin{align}
		\label{est_H_k_term}
		&\left|\dfrac{\psi_{k-2}(r)}{\tH_k(r)}\right| \leq 
			c 2^{k},\\
		\label{est_DH_k_term}
		&\left|
			\dfrac{d }{d r}
			\left(\dfrac{\psi_{k-2}(r)}{\tH_k(r)}\right)
		 \right| \leq c 2^{5k},
	\end{align}
	where $\psi_k$ are defined in \eqref{psik.unit.part.def}, $\tH_k$ is defined in \eqref{ghk_def}.
	\item[(iii)] for $k\geq 3$ and $r\geq 1-2^{-k}$ the following estimates hold:
	\begin{align}
		\label{estimate_G}
		&|\tG(r)| \leq c \dfrac{(2\sqrt{2})^{-k}}{k!}, \\
		\label{diff_estim_G}
		&\left|
			\dfrac{d \tG}{d r}(r)
		\right| \leq c\dfrac{8^{k}}{k!},
	\end{align}
	where $\tG$ is defined in \eqref{ghk_def}.
\end{itemize}
\end{lem}
\begin{lem}\label{lem_U0_part_d}
	 Let $U_0$ be defined by \eqref{U0_def}-\eqref{ghk_def}.
	 Then the following estimates 
are valid:
\begin{align}\label{U0_partialD_estim}
	&\left|
		\dfrac{\partial U_0}{\partial s}(s,r)
	\right| \leq \dfrac{C}{(1-r)^3}, \, 
	\left|
		\dfrac{\partial U_0}{\partial r}(s,r)
	\right| \leq \dfrac{C}{(1-r)^5} \text{ for }s\in \R, \, r\in (1/2,1),
\end{align}
where $C$ is a constant depending only on $\Phi$ of \eqref{phi.suppint}-\eqref{phi.unity.cond}.

\end{lem}
 Lemmas~\ref{main_lemm},~\ref{lem_U0_part_d} are proved in Subsections~\ref{subsect_main_lemm},
~\ref{subsect_lem_U0_part_d}, respectively.

\paragraph*{\large Proof of \eqref{w.w0.conv}.}
\par From \eqref{est_H_k_term}, \eqref{estimate_G} it follows that 
	\begin{align}\label{inequality_G}
		&|\tG(r)| \leq c (2\sqrt{2})^{-k+3}/(k-3)!, \\
		\label{inequality_Hk}
		&\left|\dfrac{\psi_{k-2}(r)}{\tH_k(r)}\right| \leq 
	c2^k,\\ \nonumber
		&\text{ for } r \in \Lambda_k,\, k\geq \max(4, k_1),
	\end{align}
	where $\Lambda_k$ is defined in \eqref{lambdak_def}.
	\par Properties \eqref{item2_supp.psi.1}-\eqref{item2_supp.psi.3} and estimate \eqref{est_H_k_term} imply that
	\begin{align}
	  \label{ineq_left}
	  &\begin{cases}
		\psi_{k-1}(r) = 0,\\
		\left|
			\dfrac{\psi_{k-3}(r)}{\tH_{k-1}(r)}
		\right| \leq c2^{k-1}
	   \end{cases} \text{ if } \quad r\in (1-2^{-k+3}, 1-2^{-k+2}),\\
	   &\begin{cases}\label{ineq_right}
	   		\psi_{k-2}(r)=0,\\
	   		\left|
			\dfrac{\psi_{k-1}(r)}{\tH_{k+1}(r)}
			\right| \leq c2^{k+1}
	   \end{cases}
	   \text{ if } \quad r\in (1-2^{-k+2}, 1-2^{-k+1}),\\
	   &\begin{cases}\label{ineq_mid}
	   		\psi_{k-1}(r) = 0,\\
	   		\psi_{k-3}(r)=0\\
	   \end{cases}
	   \quad \quad \quad \, 
	   \text{ if } \quad r=1-2^{-k+2},
	\end{align}
	for $k\geq \max(4,k_1)$.
	\par Note that the assumption that $r\in \Lambda_k$ is splitted into 
	the assumptions on $r$ of \eqref{ineq_left}, \eqref{ineq_right}, \eqref{ineq_mid}.
	\par Using formulas \eqref{item2_U0_long}, \eqref{inequality_G}-\eqref{ineq_mid},
	we obtain the following estimates:
	\begin{align}\label{big.estim.left}\begin{split}
		|1-U_0(s,r)| &= |\tG(r)|
		\left|
			(k-1)!\tf_{k-1}((s^2 + r^2)^{1/2})\dfrac{\psi_{k-3}(r)}{\tH_{k-1}(r)} + 
			k!\tf_{k}((s^2+r^2)^{1/2})\dfrac{\psi_{k-2}(r)}{\tH_{k}(r)}
		\right|\\
		&\leq c (2\sqrt{2})^{-k+3}(c(k-2)(k-1) 2^{k-1} + c
		(k-2)(k-1)k	2^{k})\\
		&\leq 2^5\sqrt{2}c^22^{-k/2}k^3 \quad \text{ if } \quad
		r\in (1-2^{-k+3}, 1-2^{-k+2}),
	 	\end{split}
	\end{align}
	\begin{align}\label{big.estim.right}
		\begin{split}
		|1-U_0(s,r)| &= |\tG(r)|
		\left|
			k!\tf_{k}((s^2 + r^2)^{1/2})\dfrac{\psi_{k-2}(r)}{\tH_{k}(r)}
			+(k+1)!\tf_{k+1}((p^2 + r^2)^{1/2})\dfrac{\psi_{k-1}(r)}
			{\tH_{k+1}(r)}
		\right|\\
		&\leq c(2\sqrt{2})^{-k+3}(c2^k(k-2)(k-1)k 
		+ c2^{k+1}(k-2)(k-1)k(k+1))\\
		&\leq 2^{10}\sqrt{2}c^22^{-k/2}k^4 \quad \text{ if } \quad
		r\in (1-2^{-k+2}, 1-2^{-k+1}),
		\end{split}
	\end{align}
	\begin{align}\label{big.estim.mid}
	\begin{split}
		\hspace*{-2.0cm}|1-U_0(s,r)| &= |\tG(r)|
		\left|
			k!\tf_{k}((s^2 + r^2)^{1/2})\dfrac{\psi_{k-2}(r)}{\tH_{k}(r)}
		\right|\hspace{4.2cm}\\
		&\leq 2^4\sqrt{2}c^22^{-k/2}k^3	\quad \text{ if } \quad
		r=1-2^{-k+2},
	\end{split}
	\end{align}
	for  $s\in \R, \, k\geq \max(4,k_1)$.
	Estimates \eqref{big.estim.left}-\eqref{big.estim.mid} imply that 
	\begin{align}\label{discr.estim}
		|1-U_0(s,r)|\leq C\,  2^{-k/2}k^4, \, r\in \Lambda_k, \, s\in \R, \, 
		k\geq \max(4,k_1),
	\end{align}
	where $C$ is a positive constant depending on $c$ of
	Lemma~\ref{main_lemm}. 
	\par In addition, for $r\in \Lambda_k$ we have that $2^{-k+1}< (1-r) < 2^{-k+3}$, which together with \eqref{discr.estim} imply \eqref{w.w0.conv}.
	\par This completes the proof of \eqref{w.w0.conv}.
\paragraph*{\large Proof of \eqref{U0_holder_estimate}.} We consider the following cases 
of $s,\,  s', \, r, \, r'$ in \eqref{U0_holder_estimate}:
\begin{enumerate}
	\item Let 
	\begin{equation} \label{assumption_u0_arg_c1}
		s, s'\in \R \text{ and } r,r'\geq 1. 
	\end{equation}
	Due to \eqref{U0_def} we have 
	that 
	\begin{equation}\label{U0_case_rbig}
		U_0(s,r) = 1, \, U_0(s',r') = 1.
	\end{equation}
	Identities in \eqref{U0_case_rbig} and assumption \eqref{assumption_u0_arg_c1} imply 
	\eqref{U0_holder_estimate} for this case.
	\item Let
	\begin{equation} \label{assumption_u0_arg_c2}
		s,s'\in \R, 1/2 < r < 1 \text{ and } r' \geq 1. 
	\end{equation}
		Then, due to \eqref{U0_def}, \eqref{w.w0.conv} we have that 
	\begin{align}\label{U0_case_mid_lestim}
		&|1-U_0(s,r)| \leq C(1-r)^{1/3},\\
		\label{U0_case_mid_restim}
		&U_0(s',r') = 1, 
	\end{align}
	where $s,s',r,r'$ satisfy assumption \eqref{assumption_u0_arg_c2}, $C$ is a constant depending only on $\Phi$. In particular, inequality \eqref{U0_case_mid_lestim} follows from \eqref{w.w0.conv} due to the following simple property of the logarithm: 
	\begin{equation}\label{log_property}
		\log_2^a\left(
			\frac{1}{1-r}
		\right) \leq C(a,\varepsilon) (1-r)^{-\varepsilon} \text{ for any }
		\varepsilon > 0, \, r\in [0,1), \, a > 0,
	\end{equation}
	where $C(a,\varepsilon)$ is some positive constant depending only on $a$ and $\varepsilon$.
	\par Due to \eqref{assumption_u0_arg_c2}, \eqref{U0_case_mid_lestim}, \eqref{U0_case_mid_restim} 	we have that 
	\begin{align}\label{U0_case_mid_festim}
	\begin{split}
		|U_0(s',r') - U_0(s,r)| &= |1 - U_0(s,r)| \leq C(1-r)^{1/3} \\
		&\leq C|r-r'|^{1/3} \leq C(|r-r'|^{1/3} + |s-s'|^{1/3}),
	\end{split}
	\end{align}
	where $C$ is a constant depending only on $\Phi$.
	\par Estimate \eqref{U0_case_mid_festim} and assumptions \eqref{assumption_u0_arg_c2} 
	imply \eqref{U0_holder_estimate} for this case.

	\item Let 
		\begin{equation}\label{assumption_u0_arg_c3}
			s,s'\in \R \text{ and } r,r'\in (1/2,1). 
		\end{equation}
		In addition, without loss of generality we  assume that $r > r'$. 
	\par Next, using \eqref{U0_smooth_property} one can see that 
	\begin{align}\label{U0_holder_init}
	\begin{split}
		|U_0(s,r) - U_0(s',r')| &= |U_0(s,r) - U_0(s',r) + U_0(s',r) - U_0(s',r')|\\
			&\leq |U_0(s,r) - U_0(s',r)| + |U_0(s',r) - U_0(s',r')| \\ 
			&\leq \left| \dfrac{\partial U_0}{\partial s}(\hat{s},r)
			\right| |s-s'| + 
			\left| 
				\dfrac{\partial U_0}{\partial r}(s',\hat{r})
			\right||r-r'|,
	\end{split}\\ \nonumber
			\text{for } s, s'\in \R, \, r, r' >& 1/2, 
			\text{ and for appropriate } \hat{s}, \, \hat{r}.
	\end{align}
	Note that $\hat{s}, \hat{r}$ belong to open intervals $(s,s'), (r',r)$, respectively. 
	\par Using \eqref{w.w0.conv}, \eqref{U0_partialD_estim}, \eqref{U0_case_mid_lestim},  \eqref{U0_holder_init} and the property that $ 1/2 < r' < \hat{r} < r < 1$ we obtain
	\begin{align}\label{U0_holder_power_estim}
		&|U_0(s,r) - U_0(s',r')| \leq C((1-r)^{1/3} + (1-r')^{1/3}),\\
		\label{U0_holder_power_estim2}
		&|U_0(s,r) - U_0(s',r')| \leq \dfrac{C}{(1-r)^5}(|s-s'| + |r-r'|),
	\end{align}
	where $C$ is a constant depending only on $\Phi$.
	\par We have that
	\begin{align}\nonumber
		(1-r)^{1/3} + (1-r')^{1/3} &= (1-r)^{1/3} + ((1-r) + (r-r'))^{1/3} \\ \nonumber
		&\leq 2(1-r)^{1/3} + |r-r'|^{1/3}\\
		&\leq \begin{cases}\label{holder_ineq_general}
			& \hspace{-0.3cm} 3|r-r'|^{1/3} \text{ if } 1-r \leq |r-r'|,\\
			& \hspace{-0.3cm} 3(1-r)^{1/3} \text{ if } 1-r > |r-r'|,
		 \end{cases}
	\end{align}
	where $r, r'$ satisfy \eqref{assumption_u0_arg_c3}.
	Note that in \eqref{holder_ineq_general} we used 
	the following inequality:
	\begin{equation}\label{sq_root_inequality}
		(a + b)^{1/m} \leq a^{1/m} + b^{1/m} \text{ for } a\geq 0, \, b\geq 0, \, m\in \mathbb{N}.
	\end{equation}
	In particular, using \eqref{U0_holder_power_estim}, \eqref{holder_ineq_general} we have that 
	\begin{align}\label{high_power_U0diff}
		|U_0(s,r)-U_0(s',r')|^{15} \leq 3^{15}C^{15}(1-r)^5 
		\text{ if } 1-r > |r-r'|,
	\end{align}
	where $s,s',r,r'$ satisfy assumption \eqref{assumption_u0_arg_c3}, 
	$C$ is a constant of \eqref{U0_holder_power_estim}, \eqref{U0_holder_power_estim2}.
	\par Multiplying the left and the right hand-sides of \eqref{U0_holder_power_estim2}, \eqref{high_power_U0diff} we obtain
	\begin{equation}\label{highest_power_U0diff}
		|U_0(s,r)-U_0(s',r')|^{16} \leq 3^{15}C^{16}(|s-s'| + |r-r'|), 
		\text{ if } 1-r > |r-r'|.
	\end{equation}
	\par Using \eqref{U0_holder_power_estim}, \eqref{holder_ineq_general}
	we obtain	 
		\begin{align}\label{holder_estimate_final1}
		|U_0(s,r) - U_0(s',r')| \leq 
			3C|r-r'|^{1/3}, \text{ if } 1-r \leq |r-r'|,
		\end{align}
		where $C$ is a constant of \eqref{U0_holder_power_estim}, \eqref{U0_holder_power_estim2} depending only on $\Phi$.
	Using \eqref{highest_power_U0diff} and \eqref{sq_root_inequality} for $m=16, \, a = |s-s'|, \, b = |r-r'|$, we have that
	\begin{align}\label{holder_estimate_final2}
	|U_0(s,r)-U_0(s',r')|\leq 3C(|s-s'|^{1/16} + |r-r'|^{1/16}), \text{ if }
	1-r > |r-r'|,
\end{align}
where $s,s',r,r'$ satisfy assumption \eqref{assumption_u0_arg_c3}, $C$ is a constant of  \eqref{U0_holder_power_estim}, \eqref{U0_holder_power_estim2} which depends only on $\Phi$.
	\par Formulas \eqref{holder_estimate_final1}, \eqref{holder_estimate_final2} imply \eqref{U0_holder_estimate} for this 
	case.
\end{enumerate}
Note that assumptions \eqref{assumption_u0_arg_c1}, \eqref{assumption_u0_arg_c2}, \eqref{assumption_u0_arg_c3} for cases 1, 2, 3, respectively, cover all possible choices of $s,s',r,r'$ in \eqref{U0_holder_estimate}. 
\par This completes the proof of \eqref{U0_holder_estimate}.
\par This completes the proof of Lemma~\ref{b.cont.lem}.

\section{Proofs of Lemmas~\ref{main_lemm},~\ref{lem_U0_part_d}}\label{pr.lm5}
\subsection{Proof of Lemma~\ref{main_lemm}}\label{subsect_main_lemm}
\paragraph{\large Proof of \eqref{estimate_func}, \eqref{estim_Dfunc}.} \par Estimates \eqref{estimate_func}, \eqref{estim_Dfunc}
 follow directly from \eqref{m.def.fk}-\eqref{phi.unity.cond}.

\paragraph{\large Proof of \eqref{estimate_G}.}
\par We will use the following parametrization of the points $y$ on $\gamma(x,\theta)\in T\Sp^1, 
\, (x,\theta)\in \Omega, \, r(x,\theta)\neq 0$ (see notations \eqref{pr.ts2}, \eqref{distance.line.func},  \eqref{ray_def_omega} for $d=2$):
\begin{align}\label{y_param}
	&y(\beta) = x-(x\theta)\theta + \tan(\beta)r(x,\theta)\,\theta, \, \beta \in (-\pi/2, \pi/2),
\end{align}
where $\beta$ is the parameter.
\par We have that:
\begin{align}\label{plain_mes}
	d\sigma (\beta) = r\,  d(\tan(\beta)) = \dfrac{r\, d\beta}{\cos^2\beta}, \, 
	r = r(x,\theta),
\end{align}
where $\sigma$ is the standard Lebesgue measure on $\gamma(x,\theta)$.

\par From definitions \eqref{m.def.ser.}, \eqref{ghk_def} it follows that
\begin{align}\label{G_def}
	&\tG(r) = \sum\limits_{k=1}^{\infty} \dfrac{\tG_k(r)}{k!},\\
	&\label{Gk_def}
	\tG_k(r) = \int\limits_{\gamma_r}\tf_k(|y|)\,  dy, \,  \gamma_r \in T(r), \, r > 1/2,	
\end{align}
where $T(r)$ is defined by \eqref{t_def}.
\par Using \eqref{m.def.fk}, \eqref{y_param}, \eqref{plain_mes}, \eqref{Gk_def} we obtain the following formula for $\tG_k$:
\begin{align}\nonumber
	\tG_k(r) &= r\int\limits_{-\pi/2}^{\pi/2} 
	\Phi\left(2^{k}\left(1-\frac{r}{\cos\beta}\right)\right) \cos \left(
		8^k \dfrac{r^2}{\cos^2\beta}
	\right)\dfrac{d\beta}{\cos^2\beta}\\ \nonumber
	&= \{u = \tan(\beta)\} = 2\,r\int\limits_{0}^{+\infty} 
	\Phi\left(2^{k}\left(1-r\sqrt{u^2 + 1}\right)\right) \cos \left(
		8^k r^2(u^2 + 1)
	\right)\, du\\ \nonumber
	&=\{t=u^2\}=  r \int\limits_{0}^{+\infty}
	\Phi\left(2^{k}\left(1-r\sqrt{t+1}\right)\right) \cos \left(
		8^k r^2(t+1)
	\right)\dfrac{dt}{\sqrt{t}}\\ \nonumber
	&= r\cos(8^kr^2)
	\int\limits_{0}^{+\infty}\Phi(2^k(1-r\sqrt{t+1}))
	\dfrac{\cos(8^kr^2t)}{\sqrt{t}}dt\\ \nonumber 
	&- r\sin(8^kr^2)\int\limits_{0}^{+\infty}
	\Phi(2^k(1-r\sqrt{t+1}))\dfrac{\sin(8^kr^2t)}{\sqrt{t}}dt\\
	\label{G_k_value}
	\begin{split}
		&=8^{-k/2}r^{-1}\cos(8^kr^2)
		\int\limits_{0}^{+\infty}\Phi_k(t, r)\dfrac{\cos(t)}{\sqrt{t}}dt\\ 
		&-8^{-k/2}r^{-1}\sin(8^kr^2)
		\int\limits_{0}^{+\infty}\Phi_k(t, r)\dfrac{\sin(t)}{\sqrt{t}}dt, 
		\, r > 1/2,
	\end{split}
\end{align}
where
\begin{align}\label{Phi_k_def}
	\begin{split}
	&\Phi_k(t,r) = \Phi(2^k(1-r\sqrt{8^{-k}r^{-2}t+1})), \,
	t\geq 0, \, r>1/2,  \, k\in \mathbb{N}.
	\end{split}
\end{align}
\par For integrals arising in \eqref{G_k_value} the following estimates hold:
\begin{align}
	\label{estim_phi_k_sin}
	&\left|
		\int\limits_{0}^{+\infty}\Phi_k(t,r)\dfrac{\sin(t)}{\sqrt{t}}dt
	\right| \leq C_1 < +\infty,\\	
	\label{estim_phi_k_cos}
	&\left|
		\int\limits_{0}^{+\infty}\Phi_k(t,r)\dfrac{\cos(t)}{\sqrt{t}}dt
	\right| \leq C_2 < +\infty, \\ \nonumber
	&\text{for } 1/2 < r < 1, \, k \geq 1.
\end{align}
where $\Phi_k$ is defined in \eqref{Phi_k_def}, $C_1,C_2$ are some positive constants
depending only on $\Phi$ and not depending on $k$ and $r$.
\par Estimates \eqref{estim_phi_k_sin}, \eqref{estim_phi_k_cos} are proved in 
Subsection~\ref{sect_estimates_ksincos}.
\par From \eqref{G_k_value}-\eqref{estim_phi_k_cos} it follows that
\begin{equation}\label{G_k_uniform_bound}
	|\tG_k(r)| \leq 2\cdot 8^{-k/2}(C_1 + C_2) \text{ for } r > 1/2, \, k\in \mathbb{N}.
\end{equation} 
\par  Note that for $y\in \gamma_r$, the following inequality holds:
\begin{equation}\label{cross.estim}
	\begin{split}
	&2^k(1-|y|) \leq 2^k(1 - r) \leq 2^{k-m} \leq 1/2 < 4/5\\
	&\text{for } 1-2^{-m} \leq r  < 1, \, k < m, \, m\geq 3,
	\end{split}
\end{equation}
where $\gamma_r$ is a ray in $T(r)$ (see notations of \eqref{t_def}, $d=2$).
\par Formulas \eqref{m.def.fk}, \eqref{phi.suppint}, 
\eqref{ghk_def}, \eqref{cross.estim} imply that
\begin{equation}\label{G_k_zero}
	\gamma_r\cap \supp\, f_k = \emptyset \text{ if } 
	r \geq 1-2^{-m},\,  k < m,
\end{equation}

\par In turn, \eqref{Gk_def}, \eqref{G_k_zero} imply that 
\begin{equation}\label{G_k_prop2}
	\tG_k(r) = 0 \text{ for } r \geq 1-2^{-m}, \, k < m, \, 
	m\geq 3.
\end{equation}
\par Due to \eqref{G_def}, \eqref{Gk_def}, \eqref{G_k_uniform_bound}, \eqref{G_k_prop2} we have that:
\begin{align}
	\begin{split}
	|\tG(r)| &\leq  \sum\limits_{k=m}^{\infty} |\tG_k(r)|/k! \\
	&\leq 2(C_1 + C_2)\dfrac{(2\sqrt{2})^{-m}}{m!} \sum\limits_{k=0}^{\infty} (2\sqrt{2})^{-k}
	= c_1 \dfrac{(2\sqrt{2})^{-m}}{m!}, \\
	c_1 = &(C_1+C_2)\dfrac{4\sqrt{2}}
	{2\sqrt{2}-1}, 
	\end{split}\\ \nonumber
	\text{for }r &\geq 1-2^{-m}, \, m\geq 3.
\end{align} 
This completes the proof of estimate \eqref{estimate_G}.
\paragraph{\large Proof of \eqref{diff_estim_G}.}
\par Using \eqref{G_def}, \eqref{Gk_def} we have that:
\begin{align}\label{Gk_der_series}
	\left|
		\dfrac{d \tG}{d r} (r)
	\right| \leq \sum\limits_{k=1}^{\infty}\dfrac{1}{k!}
	\left|
		\dfrac{d \tG_k(r)}{d r}
	\right|.
\end{align}
\par Formulas \eqref{m.def.fk}, \eqref{U0_r_high_deriv} for $n=1$, \eqref{estim_Dfunc}, \eqref{Gk_def} imply that
\begin{align}\label{Gk_der_proof}
	\begin{split}
	\left|
		\dfrac{d \tG_k}{d r}(r)
	\right| &= \left|
		\int\limits_{-\infty}^{+\infty}
			\dfrac{r\tf_k'((s^2+r^2)^{1/2})}{\sqrt{r^2 + s^2}}\, ds 
			\right|
			 \\
	&\leq \int\limits_{-\infty}^{+\infty} 
	|\tf_k'((s^2 + r^2)^{1/2})|\, ds 
	= 
		\int\limits_{\gamma_r}
		|\tf_k'(|y|)|\, dy
	 \leq c8^k
		\hspace*{-0.4cm}
		\int\limits_{\gamma_r\cap B(0,1)}\hspace*{-0.4cm} dy \leq 2c8^k,
	\end{split}
\end{align}
where $B(0,1)$ is defined in \eqref{ball_def}, $d=2$.\\
At the same time, formula \eqref{G_k_prop2} implies that
\begin{equation}\label{Gk_der_zero}
	\dfrac{d \tG_k(r)}{d r} = 0 \text{ for } r \geq 1-2^{-m}, \,
	k < m, \,  
	m\geq 3.
\end{equation}
Formulas \eqref{Gk_der_series}, \eqref{Gk_der_proof}, \eqref{Gk_der_zero} imply the 
following sequence of inequalities:
\begin{align}\label{Gk_final_ser_estim}
	\left|
		\dfrac{d \tG(r)}{d r}
	\right|
	\leq \sum\limits_{k=m}^{\infty}\dfrac{1}{k!}\left|
		\dfrac{d \tG_k(r)}{d r}
		\right|
	\leq c\dfrac{8^m}{m!}\sum\limits_{k=0}^{\infty}\dfrac{m!8^k}{(k+m)!}, \, 
	r \geq 1-2^{-m}, \, m\geq 3.
\end{align}
The series in the right hand-side in \eqref{Gk_final_ser_estim} admits 
the following estimate:
\begin{equation}\label{stirling_estimate}
	\sum\limits_{k=0}^{\infty}\dfrac{m!8^k}{(k+m)!} \leq \sum\limits_{k=0}^{\infty}
	\dfrac{8^k}{k!} = e^{8} 
	\text{ and the estimate does not depend on } m.
\end{equation}
\par Formulas \eqref{Gk_final_ser_estim}, \eqref{stirling_estimate} imply \eqref{diff_estim_G}.
 
\paragraph{\large Proof of \eqref{est_H_k_term}.}
\par  For each $\psi_k$ from \eqref{psik.unit.part.def} we have that
\begin{equation}
	|\psi_k| \leq 1.
\end{equation}
Therefore, it is sufficient to show that 
\begin{equation}
	\tH_k \geq C2^{-k} \text{ for } k\geq k_1, \, C = c^{-1}.
\end{equation}
Proceeding from \eqref{ghk_def} and in a similar way with \eqref{G_k_value} we obtain 
the formulas
\begin{align} \label{Hk_def_integr}
	&\tH_k(r) = r\int\limits_{0}^{+\infty}
	\dfrac{\Phi^2(2^k(1-r\sqrt{t+1}))}{\sqrt{t}} \cos^2(8^kr^2(t+1))\, dt = \tH_{k,1}(r) + \tH_{k,2}(r), \, 
	r > 1/2,
	\\ 
	\label{Hk1_def}
	&\tH_{k,1}(r) = \dfrac{r}{2}\int\limits_{0}^{+\infty} 
	\dfrac{\Phi^2(2^k(1-r\sqrt{t+1}))}{\sqrt{t}}\, dt, \\
	\label{Hk2_def} 
	&\tH_{k,2}(r) = \dfrac{r}{2}\int\limits_{0}^{+\infty}
	\dfrac{\Phi^2(2^k(1-r\sqrt{t+1}))}{\sqrt{t}} \cos(2\cdot 8^kr^2(t+1)) \, dt.
\end{align}
In addition, we have that:
\begin{equation}\label{Hk_supp.condition}
	\supp_t \Phi^2(2^{k}(1-r\sqrt{t+1}) \subset [0,3] \text{ for } 1/2 < r \leq 1-2^{-k+1}, \, 
	k\geq 3,
\end{equation}
where $\supp_t$ denotes the support of the function in variable $t$. Property \eqref{Hk_supp.condition} is proved below in this paragraph (see formulas \eqref{phi2.positive.supp}-\eqref{u1u2_upper_bound}).

\par Note that 
\begin{align}\label{secnd.assump}
	&2^k(1-r) \geq 2^{k}\cdot 2^{-k+1} \geq 2 > 6/5 \text{ for }
	1/2 < r\leq 1-2^{-k+1}, \, k \geq 3.
\end{align}
From \eqref{phi.suppint}, \eqref{phi.supp} and from \eqref{secnd.assump} we have that:
\begin{equation}\label{phi2.positive.supp}
	\supp_t \Phi^2(2^{k}(1-r\sqrt{t+1}) \subset [0,+\infty) \text{ for } 1/2 < r \leq 1-2^{-k+1}, \, 
	k\geq 3.
\end{equation}

\noindent We have that 
\begin{align}\label{u1u2_exist}
	\begin{split}
	&\exists t^{(k)}_1 = t^{(k)}_1(r)\geq 0, t^{(k)}_2=t^{(k)}_2(r)\geq 0, \,  t^{(k)}_2 > t^{(k)}_1,  \text{ such that } \\
	&\begin{cases}
		2^k(1-r\sqrt{t^{(k)}_1+1}) = 11/10,\\
		2^k(1-r\sqrt{t^{(k)}_2+1}) = 9/10,
	 \end{cases}
	 \end{split}\\
	\label{u1u2_length}
	&|t^{(k)}_2-t^{(k)}_1| \geq \left(\sqrt{t^{(k)}_2 + 1} - \sqrt{t^{(k)}_1 + 1}\right) 
	= \frac{2^{-k}}{5}
	r^{-1} \geq \frac{2^{-k}}{5}, \\ \nonumber
	&\text{for } 1/2 < r\leq 1-2^{-k+1}, \, k \geq 3.
\end{align}
In addition, from \eqref{u1u2_exist} it follows that 
\begin{align}\label{u1u2_upper_bound}
	\begin{split}
		&t^{(k)}_1 = \dfrac{(1-2^{-k}\frac{11}{10})^2}{r^2} - 1 \leq 
		4(1-2^{-k}\frac{11}{10})^2 -1 \leq 3,\\
		&t^{(k)}_2 = \dfrac{(1-2^{-k}\frac{9}{10})^2}{r^2} - 1 \leq 
		4(1-2^{-k}\frac{11}{10})^2 -1 \leq 3, 
	\end{split}\\ \nonumber 
	&\text{for } 1/2 < r \leq 1 - 2^{-k+1}, \, k \geq 3.
\end{align}
\noindent Using \eqref{phi.suppint}-\eqref{phi.unity.cond}, \eqref{Hk1_def}, \eqref{Hk_supp.condition}, \eqref{u1u2_exist}-\eqref{u1u2_upper_bound} we have that 
\begin{align}\nonumber
	\tH_{k,1}(r) &\geq  \dfrac{r}{2}
	\int\limits_{t_1^{(k)}}^{t_2^{(k)}} \dfrac{dt}{\sqrt{t}} \geq 
	\dfrac{r}{2}\hspace{-0.6cm}
	\int\limits_{3}^{3 + |t^{(k)}_2-t^{(k)}_1|} 
	\hspace{-0.6cm} \dfrac{dt}{\sqrt{t}} \\ \label{mean_val_thm}
	&\geq \dfrac{r}{6} \hspace{-0.6cm}
	\int\limits_{3}^{3+|t_2^{(k)}-t^{(k)}_1|} \hspace{-0.6cm} dt 
	= \dfrac{r}{6}|t_2^{(k)}-t^{(k)}_1| \geq \dfrac{2^{-k}}{30}
	\text{ for }\, \, 1/2 < r \leq  1-2^{-k+1}, \, k\geq 3.
\end{align}

\noindent On the other hand, proceeding from using \eqref{Hk2_def} and, in a similar way with  \eqref{G_k_value}-\eqref{G_k_uniform_bound},   
we have
\begin{align}\label{last_bound}
	\begin{split}
	 |\tH_{k,2}(r)| &= 
	\dfrac{r}{2}
	\left|\int\limits_{0}^{+\infty}
	\dfrac{\Phi^2
	(2^k(1-r\sqrt{t+1}))}{\sqrt{t}} 
	\cos(2\cdot 8^kr^2(t+1)) \, dt
	\right| \\
	&\leq \dfrac{r}{2} |\cos(2\cdot 8^kr^2)| \left|
		\int\limits_{0}^{+\infty}\Phi^2(2^k(1-r\sqrt{t+1}))
		\dfrac{\cos(2\cdot 8^kr^2t)}{\sqrt{t}}dt
	\right| \\
	&+ \dfrac{r}{2} |\sin(2\cdot 8^kr^2)| 
	\left|
		\int\limits_{0}^{+\infty}\Phi^2(2^k(1-r\sqrt{t+1}))
		\dfrac{\sin(2\cdot 8^kr^2t)}{\sqrt{t}}dt
	\right|\\
	&\leq 8^{-k/2}\dfrac{r^{-1}}{2}
		\left|
			\int\limits_{0}^{+\infty}\Phi^2_k(t,r)\dfrac{\cos(2t)}{\sqrt{t}}dt
		\right|
	+8^{-k/2}\dfrac{r^{-1}}{2}
		\left|
			\int\limits_{0}^{+\infty}\Phi^2_k(t,r)\dfrac{\sin(2t)}{\sqrt{t}}dt
		\right|\\
	&\leq 8^{-k/2}C, \text{ for } 1/2 < r < 1-2^{-k+1}, \, k\geq 3, 
	\end{split}
\end{align}
where $\Phi_k(t,r)$ is defined in \eqref{Phi_k_def},  $C$ is some constant depending only on $\Phi$ and not depending on $k, r$. In \eqref{last_bound} we have also used that $\Phi^2(t)$ satisfies assumptions \eqref{phi.suppint}-\eqref{phi.unity.cond}.

Note also that $\Phi^2(t)$ satisfies assumptions \eqref{phi.suppint}-\eqref{phi.unity.cond} for $\Phi(t)$. 
\par Using \eqref{Hk_def_integr}-\eqref{Hk2_def}, \eqref{mean_val_thm}, \eqref{last_bound} we obtain 
\begin{align}\label{H_ineq}
\begin{split}
	|\tH_k(r)| &\geq |\tH_{k,1}(r)| - |\tH_{k,2}(r)| \\ 
	&\geq \dfrac{2^{-k}}{30} - C'\cdot 8^{-k/2}\\
	&\geq 2^{-k}\left(\dfrac{1}{30} - \dfrac{C'}{(\sqrt{2})^k}\right)\\
	&\geq C2^{-k} \, \text{ for } 1/2 < r < 1-2^{-k+1}, \, k\geq 
	k_1 \geq 3, \\
	C &= \dfrac{1}{30} - C'(\sqrt{2})^{-k_1},
\end{split}
\end{align}
where $C'$ depends only on $\Phi$, $k_1$ is arbitrary constant such that $k_1 \geq 3$ and $C$ is positive. 
\par Formulas \eqref{est_H_k_term} follows from \eqref{psik.unit.part.def}, \eqref{H_ineq}.
\par This completes the proof \eqref{est_H_k_term}.

\paragraph{\large Proof of \eqref{est_DH_k_term}.}
\par The following formula holds:
\begin{equation}\label{Hk_complex_derivative}
	\dfrac{d }{d r}\left(
		\dfrac{\psi_{k-2}(r)}{\tH_k(r)}
	\right) = -
	\dfrac{\tH_k'(r)\psi_{k-2}(r) - \tH_k(r)\psi_{k-2}'(r)}{\tH_k^2(r)}, \, 1/2 < r < 1 ,
\end{equation}
where $\tH_k',\,  \psi_{k-2}'$ denote the derivatives of $\tH_k, \, \psi_k$, defined in \eqref{ghk_def}, \eqref{psik.unit.part.def}, respectively.
\par Using \eqref{m.def.fk}, \eqref{ghk_def}, \eqref{U0_r_high_deriv}, $n=1$, \eqref{estimate_func}, \eqref{estim_Dfunc}
we have that
\begin{align}\label{Hk_der_estim}
	\begin{split}
	|\tH_k'(r)| &= 
	2\left|
		\int\limits_{-\infty}^{+\infty}
		\dfrac{r}{\sqrt{r^2+s^2}}\tf_k(\sqrt{r^2 + s^2})\tf_k'(\sqrt{r^2+s^2})\, ds
	\right| \\
	&\leq 2
		\int\limits_{-\infty}^{+\infty}
		\left|
		\tf_k(\sqrt{r^2 + s^2})\tf_k'(\sqrt{r^2+s^2})
		\right|
		\, ds
	= 2
			\int\limits_{\gamma_r}|\tf_k(|y|)\tf_k'(|y|)|\, dy
			 \\
	& \leq 2c8^{k}\hspace*{-0.4cm}
	\int\limits_{\gamma_r\cap B(0,1)} \hspace*{-0.4cm}
	dy \leq 4c8^{k}, \, \gamma_r\in T(r),\, k\geq 3, \, r> 1/2,
	\end{split}
\end{align}
where we use notations \eqref{t_def}, \eqref{ball_def}, $d=2$.
\par Using \eqref{psik.unit.part.def}, \eqref{psik.unit.part.deriv}, \eqref{est_H_k_term}, \eqref{H_ineq}-\eqref{Hk_der_estim} we have that
\begin{align}
	&\left|\dfrac{d }{d r}\left(
		\dfrac{\psi_{k-2}(r)}{\tH_k(r)}
	\right)\right| \leq C2^{2k}(|\tH_k'(r)| + |\tH_k(r)|\cdot |\psi_k'(r)|)
					\leq C' 2^{5k}, \\ \nonumber
	&\text{for } 1/2 < r < 1-2^{-k+1}, \, k\geq k_1 \geq 3,
\end{align}
where $C'$ is a constant not depending on $k$ and $r$ and depending only on $\Phi$.
\par This completes the proof of Lemma~\ref{main_lemm}.

\subsection{Proof of Lemma~\ref{lem_U0_part_d}}\label{subsect_lem_U0_part_d}
It is sufficient to show that 
\begin{align}
	&\label{estimD_U0s_estimr}
	\left|
		\dfrac{\partial U_0(s,r)}{\partial s}
	\right| \leq \dfrac{C}{(1-r)^3}, \\
	\label{estimD_U0r_estimr}
	&\left|
		\dfrac{\partial U_0(s,r)}{\partial r}
	\right| \leq \dfrac{C}{(1-r)^5}, \\ \nonumber
	&\text{for } s\in \R, \, r\in \Lambda_k, \, k\geq \max(4,k_1),
\end{align}
where $C$ is a positive constant depending only on $\Phi$ of \eqref{m.def.fk}, 
$\Lambda_k$ is defined in \eqref{lambdak_def}, $k_1$ is a constant arising in Lemma~\ref{main_lemm} and depending only on $\Phi$.
\par Indeed, estimates \eqref{U0_partialD_estim} follow from \eqref{U0_smooth_property}, \eqref{estimD_U0s_estimr}, \eqref{estimD_U0r_estimr} and the fact that ${\Lambda_k,\, k\geq 4}$, is an open cover of $(1/2,1)$.
\par In turn, estimates \eqref{estimD_U0s_estimr}, \eqref{estimD_U0r_estimr} follow from 
the estimates
\begin{align}
	&\label{estimD_U0s_estimk}
	\left|
		\dfrac{\partial U_0(s,r)}{\partial s}
	\right| \leq C \cdot 8^k,\\
	&\label{estimD_U0r_estimk}
	\left|
		\dfrac{\partial U_0(s,r)}{\partial r}
	\right| \leq C\cdot  (32)^{k},\\ \nonumber
	&\text{for } s\in \R, \, r\in \Lambda_k,
\end{align}
and from the fact that $2^{-k+1}< 1-r < 2^{-k+3}, \, k\geq \max(4, k_1)$, for  $r\in \Lambda_k$, 
where $C$ is a positive constant depending only on $\Phi$.
\par Estimate \eqref{estimD_U0s_estimk} follows from formula \eqref{U0_p_high_deriv} for $n=1$ 
and estimates \eqref{estim_Dfunc}, \eqref{est_H_k_term}, 
\eqref{inequality_G}-\eqref{ineq_mid}.
\par Estimate \eqref{estimD_U0r_estimk} follows from \eqref{item2_U0_long}, 
\eqref{estimate_func}-\eqref{est_DH_k_term}, \eqref{inequality_G}-\eqref{ineq_mid} and from 
the estimates:
\begin{align}\label{derivativeH_terms_estimates}
&\left|\dfrac{d}{dr}
\left(
	\dfrac{\psi_{k-i}(r)}{\tH_{k-i+2}(r)}
\right)\right|
 \leq c 2^{5(k+1)},\\	\label{derivativeG_global}
&\left| \dfrac{d\tG(r)}{dr} \right| \leq c\dfrac{8^{-k+3}}{(k-3)!},\\ \nonumber
&\text{for } r\in \Lambda_k, \, i\in \{1,2,3\},
\end{align}
where $c$ is a constant arising in Lemma~\ref{main_lemm}.
\par Estimate \eqref{derivativeH_terms_estimates} follows from \eqref{est_DH_k_term} (used with $k-1,\, k, \, k+1$ in place of $k$).
Estimate \eqref{derivativeG_global} follows from \eqref{diff_estim_G} (used with $k-3$ in place of $k$).
\par This completes the proof of Lemma~\ref{lem_U0_part_d}.

\subsection{Proof of estimates \eqref{estim_phi_k_sin}, \eqref{estim_phi_k_cos}}
We use the following Bonnet's integration formulas (see, e.g., \cite{fichtenholz1959intdiff}, 
Chapter 2):
\begin{align}\label{bonnet_decreasing}
	&\int\limits_{a}^{b}f_1(t)h(t) \,dt = f_1(a)
	\int\limits_{a}^{\xi_1}h(t)\,dt, \\
	\label{bonnet_increasing}
	&\int\limits_{a}^{b} f_2(t)h(t)\, dt = f_2(b)
	\int\limits_{\xi_2}^{b}h(t)\, dt,
\end{align}
for some appropriate $\xi_1, \, \xi_2\in [a,b]$, where 
\begin{align}
	\label{bonnet_properties}
	\begin{split}
	&f_1\text{ is monotonously decreasing on }[a,b], \, f_1 \geq 0, \\
	&f_2\text{ is monotonously increasing on }[a,b], \, f_2 \geq 0, \\
	&h(t)\text{ is integrable on }[a,b].  
	\end{split}
\end{align}
\label{sect_estimates_ksincos}
Let 
\begin{align}\label{g1g2_def}
	&g_1(t) = \dfrac{\sin(t)}{\sqrt{t}}, \, g_2(t)=\dfrac{cos(t)}{\sqrt{t}}, \, 
	t > 0,\\
	\label{G1G2_def}
	&G_1(s) = \int\limits_{0}^{s}\dfrac{\sin(t)}{\sqrt{t}}\,dt, \, 
	 G_2(s) = \int\limits_{0}^{s}\dfrac{\cos(t)}{\sqrt{t}}\,dt, \,  s\geq 0.
\end{align}
We recall that 
\begin{equation}\label{G1G2_integrability}
	\lim\limits_{s\rightarrow +\infty} G_1(s) = \lim\limits_{s\rightarrow +\infty} G_2(s)
	= \sqrt{\dfrac{\pi}{2}}.
\end{equation}
From \eqref{g1g2_def}, \eqref{G1G2_def}, \eqref{G1G2_integrability} it follows that 
\begin{equation}\label{G1G2_bounded}
	G_1, G_2 \text{ are continuous and bounded on } [0,+\infty).
\end{equation}
Due to \eqref{phi.suppint}-\eqref{phi.monotonicity.cond}, \eqref{Phi_k_def} and monotonicity 
of the function $2^{k}(1-r\sqrt{8^{-k}r^{-2}t + 1})$ in $t$ on $[0,+\infty)$
it follows that 
\begin{align}\label{Phi_k_decreasing}
	&\Phi_k(t,r) \text{ is monotonously decreasing on } [0, +\infty), \text{ if } 2^k(1-r) \leq 11/10,\\ \label{Phi_k_decreasing_increasing}
	\begin{split}
	&\Phi_k(t,r) \text{ is monotonously increasing on } [0, t_0] \text{ for some } t_0>0\\
	&\hspace{0.55cm}\text{and is monotonously decreasing on } [t_0,+\infty), \text{ if } 
	2^k(1-r) > 11/10.
	\end{split}\\ \nonumber
	&\text{for } r > 1/2, \, k\in \mathbb{N},
\end{align}
Moreover, due to \eqref{phi.suppint}-\eqref{phi.unity.cond}, \eqref{Phi_k_def}, for $T_k=8^k, \, k\in \mathbb{N}$, we have that
\begin{align}
	&\Phi_k(T_k,r) = 
	\Phi(2^k(1-r\sqrt{r^{-2}+1})) = \Phi(2^k(1-\sqrt{1 + r^{-2}})) = 0, \\
	\label{Phi_k_right_border}
	&\Phi_k(t,r) = 0 \text{ for } t\geq T_k, \\
	&\label{Phi_k_bound}
	|\Phi_k(t,r)| \leq 1 \text{ for } t\geq 0,\\ \nonumber
	&r > 1/2, \, k\in \mathbb{N}.
\end{align} 
Using \eqref{Phi_k_def}, \eqref{g1g2_def}-\eqref{Phi_k_decreasing_increasing}, \eqref{Phi_k_right_border} and \eqref{bonnet_decreasing}-\eqref{bonnet_properties} we obtain
\begin{align}\label{integrate_decrease_estim}
	\int\limits_{0}^{+\infty}
	\Phi_k(t,r)g_i(t)\, dt &= \int\limits_{0}^{T_k}
	\Phi_k(t,r)g_i(t)\, dt = \Phi_k(0,r)\int\limits_{0}^{\xi}g_i(t)\,dt \\
	\nonumber
	&=\Phi_k(0,r)G_i(\xi) \text{ for appropriate } \xi \in [0,T_k], 
	\text{ if } 2^k(1-r) \leq 11/10,\\ \nonumber
	\int\limits_{0}^{\infty}
	\Phi_k(t,r)g_i(t)\, dt &= \int\limits_{0}^{T_k}
	\Phi_k(t,r)g_i(t)\, dt = 
	\int\limits_{0}^{t_0} \Phi_k(t,r)g_i(t)\, dt + 
	\int\limits_{t_0}^{T_k}\Phi_k(t,r)g_i(t)\, dt\\ \nonumber
	&=\Phi_k(t_0,r)\int\limits_{\xi'}^{t_0}g_i(t)\,dt + 
	\Phi_k(t_0,r)\int\limits_{t_0}^{\xi''}g_i(t)\,dt \\
	\label{integrate_jump_estim}
	&=\Phi_k(t_0,r)(G_i(\xi'')-G_i(\xi')) \\ \nonumber
	&\text{ for appropriate } \xi' \in [0,t_0],\, \xi''\in [t_0,T_k], 
	\text{ if }2^k(1-r) >  11/10,
\end{align}
where $i=\overline{1,2}$.
\par Estimates \eqref{estim_phi_k_sin}, \eqref{estim_phi_k_cos} follow from 
\eqref{g1g2_def}, \eqref{G1G2_def}, \eqref{G1G2_bounded},  \eqref{Phi_k_bound}-\eqref{integrate_jump_estim}.

\fontsize{10pt}{12pt}\selectfont
\bibliographystyle{alpha}

\end{document}